\documentclass[letterpaper,11pt]{article} % Using article class
\usepackage{algorithm,algpseudocode} % Algorithm
\usepackage{amsmath}
\usepackage{amssymb}
\usepackage{booktabs} % Table \toprule, \midrule, \bottomrule
\usepackage[margin=1in]{geometry} % Adjust page margins
\usepackage{graphicx}
\usepackage{multirow}
\usepackage{subfig} % Subfloat
\usepackage{cite} % Group citations
\usepackage{hyperref} % Reference links
\begin{document}

% ------------------------------------------------------------------------------------------------ %
% ------------------------------------------------------------------------------------------------ %
% ---------------------------------------- Heading (START) --------------------------------------- %
% ------------------------------------------------------------------------------------------------ %
% ------------------------------------------------------------------------------------------------ %
% Title
\title{\textbf{Adaptive Mesh Refinement and Error Estimation Method for Optimal Control Using Direct Collocation}}
% Author
\author{
George V. Haman III\thanks{Ph.D. Candidate, Department of Mechanical and Aerospace Engineering. Email: georgehaman@ufl.edu.}
~and Anil V. Rao\thanks{Professor, Department of Mechanical and Aerospace Engineering. Email: anilvrao@ufl.edu.} \vspace{12pt} \\ {\em{University of Florida}} \\ {\em{Gainesville, FL 32611}}
}
% Date
\date{}
% Make title
\maketitle
% ------------------------------------------------------------------------------------------------ %
% ------------------------------------------------------------------------------------------------ %
% ----------------------------------------- Heading (END) ---------------------------------------- %
% ------------------------------------------------------------------------------------------------ %
% ------------------------------------------------------------------------------------------------ %

% ------------------------------------------------------------------------------------------------ %
% ------------------------------------------------------------------------------------------------ %
% --------------------------------------- Commands (START) --------------------------------------- %
% ------------------------------------------------------------------------------------------------ %
% ------------------------------------------------------------------------------------------------ %
% Math commands
\newcommand{\deriv}[2]{\ensuremath{\dfrac{\mathrm{d}{#1}}{\mathrm{d}{#2}}}} % First derivative
\newcommand{\dt}{\ensuremath{\textrm{d}t}} % dt
\newcommand{\dtau}{\ensuremath{\mathrm{d}\tau}} % dtau
\newcommand{\dzeta}{\ensuremath{\mathrm{d}\zeta}} % dzeta
\newcommand{\lcrc}[1]{\ensuremath{\left\lceil{#1}\right\rceil}} % Left ceiling right ceiling
\newcommand{\lfrf}[1]{\ensuremath{\left\lfloor{#1}\right\rfloor}} % Left floor right floor
\newcommand{\lprp}[1]{\ensuremath{\left({#1}\right)}} % Left parenthesis right parenthesis
\newcommand{\mb}[1]{\ensuremath{\bar{#1}}} % Math bar
\newcommand{\mbb}[1]{\ensuremath{\mathbb{#1}}} % Mathbb font
\newcommand{\mbc}[1]{\ensuremath{\check{\bar{#1}}}} % Math bar check
\newcommand{\mbf}[1]{\ensuremath{\mathbf{#1}}} % Math bold face
\newcommand{\mbfh}[1]{\ensuremath{\hat{\mathbf{#1}}}} % Math bold face hat
\newcommand{\mbfbh}[1]{\ensuremath{\hat{\bar{\mathbf{#1}}}}} % Math bold face bar hat
\newcommand{\mbfbc}[1]{\ensuremath{\check{\bar{\mathbf{#1}}}}} % % Math bold face bar check
\newcommand{\mbfc}[1]{\ensuremath{\check{\mathbf{#1}}}} % Math bold face check
\newcommand{\mbh}[1]{\ensuremath{\hat{\bar{#1}}}} % Math bar hat
\newcommand{\mbs}[1]{\ensuremath{\boldsymbol{#1}}} % Math bold symbol
\newcommand{\mc}[1]{\ensuremath{\check{#1}}} % Math check
\newcommand{\mcl}[1]{\ensuremath{\mathcal{#1}}} % Mathcal font
\newcommand{\md}[1]{\ensuremath{\dot{#1}}}   % Math dot
\newcommand{\mdd}[1]{\ensuremath{\ddot{#1}}} % Math double dot
\newcommand{\mh}[1]{\ensuremath{\hat{#1}}} % Math hat
% PattersonRao mesh refinement commands
\newcommand{\PR}{\ensuremath{ph}}
\newcommand{\PRfull}{\ensuremath{ph}-\ensuremath{(N_{\min},N_{\max})}}
\newcommand{\PRvar}[2]{\ensuremath{ph}-\ensuremath{({#1},{#2}})}
% PattersonRao-HamanRao mesh refinement commands
\newcommand{\PRHR}{\ensuremath{ph^{*}}}
\newcommand{\PRHRfull}{\ensuremath{ph^{*}}-\ensuremath{(N_{\min},N_{\max},\texttt{s})}}
\newcommand{\PRHRvar}[3]{\ensuremath{ph^{*}}-\ensuremath{({#1},{#2},{#3})}}
% LiuRao mesh refinement commands
\newcommand{\LR}{\ensuremath{hp}}
\newcommand{\LRfull}{\ensuremath{hp}}
\newcommand{\LRvar}[0]{\ensuremath{hp}}
% LiuRao-HamanRao mesh refinement commands
\newcommand{\LRHR}{\ensuremath{hp^{*}}}
\newcommand{\LRHRfull}{\ensuremath{hp^{*}}-\ensuremath{(\texttt{s})}}
\newcommand{\LRHRvar}[1]{\ensuremath{hp^{*}}-\ensuremath{({#1})}}
% LiuRao-Legendre mesh refinement commands
\newcommand{\LRL}{\ensuremath{hp\text{-Legendre}}}
\newcommand{\LRLfull}{\ensuremath{hp\text{-Legendre}}}
\newcommand{\LRLvar}[0]{\ensuremath{hp\text{-Legendre}}}
% LiuRao-Legendre-HamanRao mesh refinement commands
\newcommand{\LRLHR}{\ensuremath{hp\text{-Legendre}^{*}}}
\newcommand{\LRLHRfull}{\ensuremath{hp\text{-Legendre}^{*}}-\ensuremath{(\texttt{s})}}
\newcommand{\LRLHRvar}[1]{\ensuremath{hp\text{-Legendre}^{*}}-\ensuremath{({#1})}}
% HamanRao mesh refinement commands
\newcommand{\HR}{\ensuremath{ph\texttt{s}}}
\newcommand{\HRfull}{\ensuremath{ph\texttt{s}}-\ensuremath{(N_{\min},N_{\max},\texttt{s})}}
\newcommand{\HRvar}[3]{\ensuremath{ph\texttt{s}}-\ensuremath{({#1},{#2},{#3})}}
% Miscellaneous commands
\setlength{\arraycolsep}{1.5pt}
\newcommand{\figureheight}{14.9em} % Figure height
\newcommand{\MATLAB}{\texttt{MATLAB}} % MATLAB
\newcommand{\ode}[1]{\texttt{ode{#1}}} % MATLAB ODE solver
% Alignment commands
\newlength{\mainside}
\newlength{\leftside}
\newlength{\centerside}
\newlength{\rightside}
\newcommand*{\mainterm}{}
\newcommand*{\leftterm}{}
\newcommand*{\centerterm}{}
\newcommand*{\rightterm}{}
\newcommand*{\term}[1]{\ensuremath{\displaystyle{#1}}}
% Bibliography commands
\newcommand{\enquote}[1]{``{#1},''}
% ------------------------------------------------------------------------------------------------ %
% ------------------------------------------------------------------------------------------------ %
% ---------------------------------------- Commands (END) ---------------------------------------- %
% ------------------------------------------------------------------------------------------------ %
% ------------------------------------------------------------------------------------------------ %

% ------------------------------------------------------------------------------------------------ %
% ------------------------------------------------------------------------------------------------ %
% ------------------------------------- Environments (START) ------------------------------------- %
% ------------------------------------------------------------------------------------------------ %
% ------------------------------------------------------------------------------------------------ %
% Table footnotes
\newenvironment{tablenotes}{\list{}{%
    \setlength{\labelsep}{0pt}%
    \setlength{\labelwidth}{0pt}%
    \setlength{\leftmargin}{0pt}%
    \setlength{\rightmargin}{0pt}%
    \setlength{\topsep}{2pt}%
    \setlength{\itemsep}{0pt}%
    \setlength{\partopsep}{0pt}%
    \setlength{\listparindent}{0em}%
    \setlength{\parsep}{0pt}}%
    \item\relax%
}{\endlist}%
% ------------------------------------------------------------------------------------------------ %
% ------------------------------------------------------------------------------------------------ %
% -------------------------------------- Environments (END) -------------------------------------- %
% ------------------------------------------------------------------------------------------------ %
% ------------------------------------------------------------------------------------------------ %

% ------------------------------------------------------------------------------------------------ %
% ------------------------------------------------------------------------------------------------ %
% --------------------------------------- Abstract (START) --------------------------------------- %
% ------------------------------------------------------------------------------------------------ %
% ------------------------------------------------------------------------------------------------ %
\begin{abstract}
An adaptive mesh refinement and error estimation method for numerically solving optimal control problems is developed using Legendre-Gauss-Radau direct collocation.  In regions of the solution where the desired accuracy tolerance has not been met, the mesh is refined by either increasing the degree of the approximating polynomial in a mesh interval or dividing a mesh interval into subintervals.  In regions of the solution where the desired accuracy tolerance has been met, the mesh size may be reduced by either merging adjacent mesh intervals or decreasing the degree of the approximating polynomial in a mesh interval.  Coupled with the mesh refinement method described in this paper is a newly developed relative error estimate that is based on the differences between solutions obtained from the collocation method and those obtained by solving initial-value and terminal-value problems in each mesh interval using an interpolated control obtained from the collocation method.  Because the error estimate is based on explicit simulation, the solution obtained via collocation is in close agreement with the solution obtained via explicit simulation using the control on the final mesh, which ensures that the control is an accurate approximation of the true optimal control.  The method is demonstrated on three examples from the open literature, and the results obtained show an improvement in final mesh size when compared against previously developed mesh refinement methods.
\end{abstract}
% ------------------------------------------------------------------------------------------------ %
% ------------------------------------------------------------------------------------------------ %
% ---------------------------------------- Abstract (END) ---------------------------------------- %
% ------------------------------------------------------------------------------------------------ %
% ------------------------------------------------------------------------------------------------ %

% ------------------------------------------------------------------------------------------------ %
% ------------------------------------------------------------------------------------------------ %
% ------------------------------------- Miscellaneous (START) ------------------------------------ %
% ------------------------------------------------------------------------------------------------ %
\renewcommand{\baselinestretch}{1.5}
\normalsize\normalfont
% ------------------------------------------------------------------------------------------------ %
% -------------------------------------- Miscellaneous (END) ------------------------------------- %
% ------------------------------------------------------------------------------------------------ %
% ------------------------------------------------------------------------------------------------ %

% ------------------------------------------------------------------------------------------------ %
% ------------------------------------------------------------------------------------------------ %
% --------------------------------- SECTION: Introduction (START) -------------------------------- %
% ------------------------------------------------------------------------------------------------ %
% ------------------------------------------------------------------------------------------------ %
\section{Introduction}\label{sec:intro}

Direct collocation methods for solving optimal control problems have become increasingly popular over the past few decades.
In a direct collocation method, the time domain of the optimal control problem is partitioned into a {\em{mesh}}.
Within each mesh interval, the state is approximated using a set of basis or trial functions (typically, these basis functions are polynomials), and constraints are enforced at a specific set of points called {\em{collocation points}}.
The resulting approximation transforms the original continuous-time optimal control problem into a finite-dimensional {\em{nonlinear programming problem}} (NLP)~\cite{Betts2020}, and the NLP is solved using well-developed software (for example, SNOPT~\cite{GillSaunders2005} or IPOPT~\cite{BieglerZavala2009}).

Direct collocation methods have historically been formulated as either $h$ methods or $p$ methods.
In an $h$ method, the order of the state approximation is fixed in each interval.
Accuracy using an $h$ method is then increased by adjusting the number and placement of mesh points.
In a $p$ method, the order of the state approximation varies in each interval, where the number and placement of mesh points are fixed.
Accuracy using a $p$ method is then increased by adjusting the order of the approximation in an interval.
While both $h$ and $p$ direct collocation methods have been used extensively, both approaches have limitations, specifically when obtaining high-accuracy solutions.
An $h$ method may require an extremely large number of intervals, whereas a $p$ method may require an unreasonably high-degree polynomial approximation.
In order to significantly reduce the size of the finite-dimensional approximation and improve the computational efficiency of solving the NLP, $hp$ methods have been developed.
In an $hp$ method, the number of intervals, placement of mesh points, and degree of the state approximation in each interval can vary.
Furthermore, $hp$ methods have gained considerable attention due to their robustness and ability to outperform both $h$ and $p$ methods in terms of computational efficiency and mesh size.

In recent years, the class of $hp$ direct {\em Gaussian quadrature collocation} methods have been developed.
In an $hp$ Gaussian quadrature collocation method, the state is approximated using a basis of polynomials (typically, Lagrange or Chebyshev polynomials) whose support points include Gaussian quadrature points, and collocation is then performed at the Gaussian quadrature points~\cite{ElnagarRazzaghi1995,FahrooRoss2001,BensonRao2006,GargHuntington2010,GargRao2011a,GargRao2011b}.
The most well-developed Gaussian quadrature collocation methods employ collocation at the {\em{Legendre-Gauss}} (LG) points~\cite{BensonRao2006,GargHuntington2010,GargRao2011a}, {\em{Legendre-Gauss-Radau}} (LGR) points~\cite{GargHuntington2010,GargRao2011a,GargRao2011b}, or {\em{Legendre-Gauss-Lobatto}} (LGL) points~\cite{ElnagarRazzaghi1995,FahrooRoss2001}.
Under certain assumptions of smoothness in the solution, Gaussian quadrature collocation methods that employ collocation at the LG or LGR points converge at an exponential rate as a function of the degree of the approximating polynomial~\cite{HagerRao2016,HagerRao2017,HagerWang2018,HagerWang2019}.
Motivated by the desire to improve computational efficiency and accuracy, this research develops a novel $hp$ Gaussian quadrature collocation method for optimal control driven by a novel error estimate.

In order to develop a practical $hp$ method, an appropriate mesh must be chosen.
Such a mesh is not known a priori and must be determined through an iterative process known as {\em{mesh refinement}}, where a new mesh is generated based on an estimate of the solution error on the current mesh and the NLP is then re-solved on the new mesh.
Various $hp$ Gaussian quadrature mesh refinement methods have been developed~\cite{DarbyRao2011a,DarbyRao2011b,PattersonRao2015,LiuRao2015,LiuRao2018,MillerRao2021}.
These $hp$ methods utilize $h$ refinement in regions of the solution that exhibit either nonsmooth or rapidly changing behavior and employ $p$ refinement in regions where the solution is changing in a smooth or slow manner.
In the case of $p$ refinement, methods for adjusting the polynomial degree include using the ratio of the maximum error and desired error tolerance~\cite{DarbyRao2011a} or using the known exponential convergence rate of a Gaussian quadrature collocation method~\cite{LiuRao2015,LiuRao2018,PattersonRao2015}.
It is noted that the $hp$ mesh refinement methods in Refs.~\cite{DarbyRao2011a,DarbyRao2011b,PattersonRao2015} do not allow for mesh size reduction.
As a result, the mesh may become unnecessarily large, and the large mesh can decrease the computational efficiency of solving the NLP. 
On the other hand, the methods described in Refs.~\cite{LiuRao2015,LiuRao2018} allow for mesh size reduction, which make it possible to obtain a smaller mesh while maintaining the desired solution accuracy.

The contributions of this paper are as follows.
The first contribution of this paper is a new mesh refinement method for both increasing and decreasing the size of the mesh.
This new $p$-then-$h$ mesh refinement strategy differs from the method of Ref.~\cite{PattersonRao2015} in that the $p$ refinement step adopts a conservative version of that described in Ref.~\cite{DarbyRao2011a} to avoid adding collocation points unnecessarily.
Furthermore, different from the method of Ref.~\cite{LiuRao2015}, the method of this paper can perform both $h$ and $p$ refinement on every mesh refinement iteration.
In addition, a new $h$-then-$p$ strategy for mesh size reduction is developed in this paper that differs fundamentally from the $p$-then-$h$ reduction approach developed in Ref.~\cite{LiuRao2015}.
Specifically, Ref.~\cite{LiuRao2015} employs power series approximations and requires the polynomial degree in adjacent intervals to be the same, whereas the new method developed in this paper for merging adjacent intervals ensures that the explicit simulation (time-marching) solution obtained over adjacent intervals is in close agreement with the solution obtained via collocation in each interval.
Consequently, the mesh size reduction method developed in this paper does not require that the polynomial degree be the same in adjacent intervals, which leads to greater mesh size reduction compared with the method of Ref.~\cite{LiuRao2015}.
Moreover, the method of this paper attempts to decrease the degree of the polynomial approximation in an interval based on the maximum error, desired error tolerance, and current degree of the polynomial approximation.
It is noted for completeness that preliminary versions of the approach developed in this paper are given in Refs.~\cite{HamanRao2024a,HamanRao2024b}.

A second contribution of this paper is a newly developed error estimate of the solution on a given mesh.
Previously developed error estimates are based on the violation in the discretized differential-algebraic constraints between the collocation points~\cite{DarbyRao2011a,DarbyRao2011b,PattersonRao2015}.
Unlike these previous error estimates, the error estimate developed in this paper is obtained as follows.
First, the dynamics are simulated in forward time across the mesh interval using the state at the start of the mesh interval.
Second, the dynamics are simulated in backward time across the mesh interval using the state at the terminus of the mesh interval.
During both forward and backward explicit simulations within a mesh interval, the control obtained via collocation is interpolated using a Lagrange polynomial.
Then, two relative error estimates are obtained for the state within a mesh interval.
The first relative error estimate is the maximum difference between the Lagrange polynomial approximation of the state in the mesh interval and the state obtained via forward explicit simulation in the mesh interval.
The second relative error estimate is the maximum difference between the Lagrange polynomial approximation of the state in the mesh interval and the state obtained via backward explicit simulation in the mesh interval.
The relative error estimate utilized in the mesh refinement method is then taken to be the maximum of these forward and backward relative error estimates.
Using this approach for estimating the error ensures that the solutions obtained via collocation and explicit simulation schemes are in agreement with one other on the final mesh.
As a result, it is unnecessary to re-validate the solution via explicit simulation after the NLP is solved on the final mesh (as is often done after mesh refinement is complete in order to verify that the solution obtained from the NLP agrees with that obtained via explicit simulation).

This paper is organized as follows. Section~\ref{sec:bolza} introduces the continuous-time Bolza optimal control problem.
Section~\ref{sec:lgr} describes the discretization of the Bolza optimal control problem given in Section \ref{sec:bolza} using the Legendre-Gauss-Radau collocation method.
Section~\ref{sec:mesh} describes the new relative error estimate along with the new mesh refinement method.
Section~\ref{sec:examples} demonstrates the method developed in Section~\ref{sec:mesh} on three examples taken from the open literature, where the key features of the new method are highlighted.
In addition, Section \ref{sec:examples} provides a comparison of the performance of the method developed in this paper against previously developed mesh refinement methods.
Finally, Section~\ref{sec:conclusion} provides conclusions on this research.

% ------------------------------------------------------------------------------------------------ %
% ------------------------------------------------------------------------------------------------ %
% ---------------------------------- SECTION: Introduction (END) --------------------------------- %
% ------------------------------------------------------------------------------------------------ %
% ------------------------------------------------------------------------------------------------ %

% ------------------------------------------------------------------------------------------------ %
% ------------------------------------------------------------------------------------------------ %
% ------------------------ SECTION: Bolza Optimal Control Problem (START) ------------------------ %
% ------------------------------------------------------------------------------------------------ %
% ------------------------------------------------------------------------------------------------ %
\section{Bolza Optimal Control Problem}\label{sec:bolza}

Without loss of generality, consider the following single-interval, continuous-time Bolza optimal control problem defined on the interval $\tau\in[-1,+1]$.
Determine the state, $\mbf{x}(\tau)\in\mbb{R}^{n_x}$, the control, $\mbf{u}(\tau)\in\mbb{R}^{n_u}$, the initial time, $t_0\in\mbb{R}$, and the terminal time, $t_f\in\mbb{R}$, that minimize the objective functional
\begin{equation}\label{eq:bolza-tau-obj}
    \mcl{J}=\mcl{M}(\mbf{x}(-1),t_0,\mbf{x}(+1),t_f)+\alpha\int_{-1}^{+1}\mcl{L}(\mbf{x}(\tau),\mbf{u}(\tau),\tau)\dtau,
\end{equation}
subject to the dynamic constraints
\begin{equation}\label{eq:bolza-tau-dyn}
    \deriv{\mbf{x}(\tau)}{\tau}-\alpha\mbf{a}(\mbf{x}(\tau),\mbf{u}(\tau),\tau)=\mbf{0},
\end{equation}
the boundary conditions
\begin{equation}\label{eq:bolza-tau-bnd}
    \mbf{b}(\mbf{x}(-1),t_0,\mbf{x}(+1),t_f)\leq\mbf{0},
\end{equation}
and the inequality path constraints
\begin{equation}\label{eq:bolza-tau-pth}
    \mbf{c}(\mbf{x}(\tau),\mbf{u}(\tau),\tau)\leq\mbf{0},
\end{equation}
where $\alpha\equiv\dt/\dtau=(t_f-t_0)/2$ is the {\em{time interval scaling factor}}.
The functions $\mcl{M}$, $\mcl{L}$, $\mbf{a}$, $\mbf{b}$, and $\mbf{c}$ are defined by the mappings
\begin{equation*}
    \arraycolsep=3.0pt
    \begin{array}{cclcl}
        \mcl{M} & : & \mbb{R}^{n_x}\times\mbb{R}\times\mbb{R}^{n_x}\times\mbb{R} & \rightarrow & \mbb{R}, \\
        \mcl{L} & : & \mbb{R}^{n_x}\times\mbb{R}^{n_u}\times\mbb{R}              & \rightarrow & \mbb{R}, \\
        \mbf{a} & : & \mbb{R}^{n_x}\times\mbb{R}^{n_u}\times\mbb{R}              & \rightarrow & \mbb{R}^{n_x}, \\
        \mbf{b} & : & \mbb{R}^{n_x}\times\mbb{R}\times\mbb{R}^{n_x}\times\mbb{R} & \rightarrow & \mbb{R}^{n_b}, \\
        \mbf{c} & : & \mbb{R}^{n_x}\times\mbb{R}^{n_u}\times\mbb{R}              & \rightarrow & \mbb{R}^{n_c}, \\
    \end{array}
\end{equation*}
where $n_x$, $n_u$, $n_b$, and $n_c$ denote the number of states, controls, boundary conditions, and inequality path constraints, respectively.
The affine transformations relating the original time interval, $t\in[t_0,t_f]$, and $\tau\in[-1,+1]$ are given by
\begin{equation}\label{eq:bolza-tau2t_t2tau}
    \begin{array}{rclcl}
        t    & \equiv & t(\tau,t_0,t_f) & = & \alpha\tau+\alpha_0, \\[0.5em]
        \tau & \equiv & \tau(t,t_0,t_f) & = & \dfrac{1}{\alpha}\lprp{t-\alpha_0}, \\
    \end{array}
\end{equation}
where $\alpha_0\equiv(t_f+t_0)/2$.

Suppose now that the domain $\tau\in[-1,+1]$ is partitioned into a {\em{mesh}} consisting of $K$ {\em{mesh intervals}} $\mcl{S}_k=[\tau_{k-1},\tau_k]\subseteq[-1,+1],\,(k=1,\ldots,K)$, where the {\em{mesh points}} are $-1=\tau_0<\tau_1<\cdots<\tau_{K-1}<\tau_K=+1$ and
\begin{equation}\label{eq:bolza-mesh}
    \bigcup_{k=1}^{K}\mcl{S}_k=[-1,+1],
    \qquad
    \mcl{S}_k\bigcap\mcl{S}_{k+1}=\{\tau_{k}\},
    \qquad
    (k=1,\ldots,K-1).
\end{equation}
The domain within each mesh interval is then mapped to the domain $\zeta\in[-1,+1]$, and the affine transformations relating $\tau\in[\tau_{k-1},\tau_k]$ and $\zeta\in[-1,+1]$ in each mesh interval are given by
\begin{equation}\label{eq:bolza-zeta2tau_tau2zeta}
    \begin{array}{rclcl}
        \tau  & \equiv & \tau(\zeta,\tau_{k-1},\tau_k) & = & \beta_k\zeta+\beta_{k0}, \\[0.5em]
        \zeta & \equiv & \zeta(\tau,\tau_{k-1},\tau_k) & = & \dfrac{1}{\beta_k}(\tau-\beta_{k0}), \\
    \end{array}
    \qquad
    (k=1,\ldots,K),
\end{equation}
where $\beta_k\equiv\dtau/\dzeta=(\tau_k-\tau_{k-1})/2$ is a {\em{mesh interval scaling factor}} and $\beta_{k0}\equiv(\tau_k+\tau_{k-1})/2$.
Finally, the Bolza optimal control problem defined by Eqs.~\eqref{eq:bolza-tau-obj}--\eqref{eq:bolza-tau-pth} is redefined in multiple-interval form as follows, where the notation $x^{(k)}$ denotes a quantity (in this case, $x$) in mesh interval $\mcl{S}_k$.
Determine the state, $\mbf{x}^{(k)}(\zeta)\in\mbb{R}^{n_x}$, and the control, $\mbf{u}^{(k)}(\zeta)\in\mbb{R}^{n_u}$, in mesh interval $\mcl{S}_k,\,(k=1,\ldots,K)$, the initial time, $t_0$, and the terminal time, $t_f$, that minimize the objective functional
\begin{equation}\label{eq:bolza-zeta-obj}
    \mcl{J}=\mcl{M}(\mbf{x}^{(1)}(-1),t_0,\mbf{x}^{(K)}(+1),t_f)+\alpha\sum_{k=1}^K\beta_k\int_{-1}^{+1}\mcl{L}(\mbf{x}^{(k)}(\zeta),\mbf{u}^{(k)}(\zeta),\zeta)\dzeta,
\end{equation}
subject to the dynamic constraints
\begin{equation}\label{eq:bolza-zeta-dyn}
    \deriv{\mbf{x}^{(k)}(\zeta)}{\zeta}-\alpha\beta_k\mbf{a}(\mbf{x}^{(k)}(\zeta),\mbf{u}^{(k)}(\zeta),\zeta)=\mbf{0},
    \qquad
    (k=1,\ldots,K),
\end{equation}
the boundary conditions
\begin{equation}\label{eq:bolza-zeta-bnd}
    \mbf{b}(\mbf{x}^{(1)}(-1),t_0,\mbf{x}^{(K)}(+1),t_f)\leq\mbf{0},
\end{equation}
the inequality path constraints
\begin{equation}\label{eq:bolza-zeta-pth}
    \mbf{c}(\mbf{x}^{(k)}(\zeta),\mbf{u}^{(k)}(\zeta),\zeta)\leq\mbf{0},
    \qquad
    (k=1,\ldots,K),
\end{equation}
and the continuity constraints
\begin{equation}\label{eq:bolza-zeta-cnt}
    \mbf{x}^{(k)}(+1)-\mbf{x}^{(k+1)}(-1)=\mbf{0},
    \qquad
    (k=1,\ldots,K-1).
\end{equation}
It is noted that the constraints of Eq.~\eqref{eq:bolza-zeta-cnt} ensure state continuity at every interior mesh point.

% ------------------------------------------------------------------------------------------------ %
% ------------------------------------------------------------------------------------------------ %
% ------------------------- SECTION: Bolza Optimal Control Problem (END) ------------------------- %
% ------------------------------------------------------------------------------------------------ %
% ------------------------------------------------------------------------------------------------ %

% ------------------------------------------------------------------------------------------------ %
% ------------------------------------------------------------------------------------------------ %
% ----------------------- SECTION: Legendre-Gauss-Radau Collocation (Start) ---------------------- %
% ------------------------------------------------------------------------------------------------ %
% ------------------------------------------------------------------------------------------------ %
\section{Legendre-Gauss-Radau Collocation}\label{sec:lgr}

The multiple-interval formulation of the continuous-time Bolza optimal control problem defined by Eqs.~\eqref{eq:bolza-zeta-obj}--\eqref{eq:bolza-zeta-cnt} is discretized using collocation at the standard LGR points~\cite{GargHuntington2010,GargRao2011b,PattersonRao2012}.
In mesh interval $\mcl{S}_k,\,(k=1,\ldots,K)$, the state is approximated using a basis of Lagrange polynomials, $\ell_j^{(k)}(\zeta),\,(j=1,\ldots,N_k+1)$, as
\begin{equation}\label{eq:lgr-xappx}
    \mbf{x}^{(k)}(\zeta)\approx\mbf{X}^{(k)}(\zeta)=\sum_{j=1}^{N_k+1}\ell_j^{(k)}(\zeta)\mbf{X}^{(k)}_j,
    \qquad
    \ell_j^{(k)}(\zeta)=\prod_{
    \begin{subarray}{c}
        l=1\\l\neq{j} \\
    \end{subarray}
    }^{N_k+1}\dfrac{\zeta-\zeta_l^{(k)}}{\zeta_j^{(k)}-\zeta_l^{(k)}},
\end{equation}
where $\mbf{X}_j^{(k)}\in\mbb{R}^{n_x}$ denotes the state approximation at $\zeta^{(k)}_j$, $\{\zeta_1^{(k)},\ldots,\zeta_{N_k}^{(k)}\}\in[-1,+1)$ is the set of $N_k$ LGR collocation points, and $\zeta_{N_k+1}^{(k)}=+1$ is a noncollocated point.
An advantage to employing a Lagrange polynomial approximation with LGR support points is that the classic Runge phenomenon is eliminated~\cite{Huntington2007}.
Differentiating $\mbf{X}^{(k)}(\zeta)$ in Eq.~\eqref{eq:lgr-xappx} with respect to $\zeta$ gives
\begin{equation}\label{eq:lgr-dxappx}
    \deriv{\mbf{x}^{(k)}(\zeta)}{\zeta}\approx\deriv{\mbf{X}^{(k)}(\zeta)}{\zeta}=\sum_{j=1}^{N_k+1}\deriv{\ell_j^{(k)}(\zeta)}{\zeta}\mbf{X}^{(k)}_j,
    \qquad
    (k=1,\ldots,K).
\end{equation}
Within each mesh interval, evaluating $\text{d}\ell^{(k)}_j(\zeta)/\dzeta$ in Eq.~\eqref{eq:lgr-dxappx} at the $N_k$ LGR points gives the elements of the {\em{LGR differentiation matrix}}, $\mbf{D}^{(k)}\in\mbb{R}^{N_k\times(N_k+1)}$, as
\begin{equation}\label{eq:lgr-dmat}
    D_{ij}^{(k)}\equiv\deriv{\ell_j^{(k)}(\zeta^{(k)}_i)}{\zeta},
    \qquad
    \left(
    \begin{array}{ll}
        i & =1,\ldots,N_k \\
        j & =1,\ldots,N_k+1 \\
        k & =1,\ldots,K \\
    \end{array}
    \right)\!.
\end{equation}
The discretization of the multiple-interval, continuous-time Bolza optimal control problem defined by Eqs.~\eqref{eq:bolza-zeta-obj}--\eqref{eq:bolza-zeta-cnt} is given as follows.
First, the objective functional of Eq.~\eqref{eq:bolza-zeta-obj} is approximated using a multiple-interval LGR quadrature.
Next, the dynamic constraints and inequality path constraints of Eqs.~\eqref{eq:bolza-zeta-dyn} and~\eqref{eq:bolza-zeta-pth}, respectively, are collocated at the LGR points in each mesh interval.
Then, the boundary conditions of Eq.~\eqref{eq:bolza-zeta-bnd} are approximated at the boundary points, while the continuity constraints of Eq.~\eqref{eq:bolza-zeta-cnt} are enforced using the approximation of the state at the terminus of an interior mesh interval and the start of the next interior mesh interval.
The resulting NLP is formally stated as follows.
Determine the state, $\mbf{X}^{(k)}_j,\,(j=1,\ldots,N_k+1)$, and the control, $\mbf{U}^{(k)}_j,\,(j=1,\ldots,N_k)$, in mesh interval $\mcl{S}_k,\,(k=1,\ldots,K)$, the initial time, $t_0$, and the terminal time, $t_f$, that minimize the objective function
\begin{equation}\label{eq:lgr-obj}
    \mcl{J}=\mcl{M}(\mbf{X}^{(1)}_1,t_0,\mbf{X}^{(K)}_{N_K+1},t_f)+\alpha\sum_{k=1}^{K}\beta_k\sum_{i=1}^{N_k}w_i^{(k)}\mcl{L}(\mbf{X}^{(k)}_i,\mbf{U}^{(k)}_i,\zeta_i^{(k)}),
\end{equation}
subject to the discretized dynamic constraints
\begin{equation}\label{eq:lgr-dyn}
    \sum_{j=1}^{N_k+1}D_{ij}^{(k)}\mbf{X}_j^{(k)}-\alpha\beta_k\mbf{a}(\mbf{X}_i^{(k)},\mbf{U}_i^{(k)},\zeta^{(k)}_i)=\mbf{0},
    \qquad
    \left(
    \begin{array}{ll}
        i & =1,\ldots,N_k \\
        k & =1,\ldots,K \\
    \end{array}
    \right)\!,
\end{equation}
the discretized boundary conditions
\begin{equation}\label{eq:lgr-bnd}
    \mbf{b}(\mbf{X}^{(1)}_1,t_0,\mbf{X}^{(K)}_{N_K+1},t_f)\leq\mbf{0},
\end{equation}
the discretized inequality path constraints
\begin{equation}\label{eq:lgr-pth}
    \mbf{c}(\mbf{X}^{(k)}_i,\mbf{U}^{(k)}_i,\zeta^{(k)}_i)\leq\mbf{0},
    \qquad
    \left(
    \begin{array}{ll}
        i & =1,\ldots,N_k \\
        k & =1,\ldots,K \\
    \end{array}
    \right)\!,
\end{equation}
and the discretized continuity constraints
\begin{equation}\label{eq:lgr-cnt}
    \mbf{X}^{(k)}_{N_k+1}-\mbf{X}^{(k+1)}_1=\mbf{0},
    \qquad
    (k=1,\ldots,K-1),
\end{equation}
where $\mbf{U}_i^{(k)}\in\mbb{R}^{n_u}$ and $w_i^{(k)}\in\mbb{R}$ are the control approximation at $\zeta_i^{(k)}$ and $i^{\text{th}}$ LGR quadrature weight, respectively, in mesh interval $\mcl{S}_k$.
It is noted that the constraints of Eq.~\eqref{eq:lgr-cnt} are implicitly satisfied by employing the same NLP decision variables for $\mbf{X}^{(k)}_{N_k+1}$ and $\mbf{X}^{(k+1)}_1,\,(k=1,\ldots,K-1)$.

% ------------------------------------------------------------------------------------------------ %
% ------------------------------------------------------------------------------------------------ %
% ------------------------ SECTION: Legendre-Gauss-Radau Collocation (END) ----------------------- %
% ------------------------------------------------------------------------------------------------ %
% ------------------------------------------------------------------------------------------------ %

% ------------------------------------------------------------------------------------------------ %
% ------------------------------------------------------------------------------------------------ %
% ----------------------- SECTION: Adaptive Mesh Refinement Method (START) ----------------------- %
% ------------------------------------------------------------------------------------------------ %
% ------------------------------------------------------------------------------------------------ %
\section{Adaptive Mesh Refinement Method}\label{sec:mesh}

In this section, the adaptive mesh refinement method of this paper is described.
First, Section~\ref{subsec:mesh-sim} provides a brief overview of numerical methods for solving a set of ordinary differential equations (ODEs).
Second, an estimate of the relative error in the solution on a given mesh is derived in Section~\ref{subsec:mesh-err}, which guides the mesh refinement process of Section~\ref{subsec:mesh-refine}.
Next, two methods for increasing the size of the mesh are presented in Sections~\ref{subsubsec:mesh-refine-pp} and~\ref{subsubsec:mesh-refine-hp}, which include increasing the degree of the polynomial approximation in an interval and dividing an interval.
Then, two methods for decreasing the size of the mesh are presented in Sections~\ref{subsubsec:mesh-refine-hm} and~\ref{subsubsec:mesh-refine-pm}, which include merging adjacent intervals and decreasing the degree of the polynomial approximation in an interval.
Finally, Section~\ref{subsec:mesh-sum} summarizes the adaptive mesh refinement method.

% ------------------------------------------------------------------------------------------------ %
% ------------------ SUBSECTION: Numerical Simulation in a Mesh Interval (START) ----------------- %
% ------------------------------------------------------------------------------------------------ %
\subsection{Numerical Simulation in a Mesh Interval}\label{subsec:mesh-sim}

The dynamic constraints defined by Eq.~\eqref{eq:bolza-zeta-dyn} associated with the multiple-interval formulation of the continuous-time Bolza optimal control problem in Section~\ref{sec:bolza} are represented by a set of ODEs.
In order to obtain solutions to optimal control problems numerically, the use of numerical methods for solving the set of ODEs is required, which are categorized as either {\em{collocation}} (implicit simulation) or {\em{time-marching}} (explicit simulation) methods.
As demonstrated by the LGR collocation method of Section~\ref{sec:lgr}, collocation schemes divide each mesh interval domain into steps, and the solution at every step in every mesh interval is obtained simultaneously; however, the objective of this section is to sequentially obtain the solution at each step in each mesh interval using a time-marching scheme.
Regardless of the numerical method chosen, solving the set of ODEs can be posed in the form of an initial-value problem (IVP) given an initial condition or a terminal-value problem (TVP) given a terminal condition.

Suppose now that the NLP of Eqs.~\eqref{eq:lgr-obj}--\eqref{eq:lgr-cnt} is solved on a given mesh, resulting in the values of the state, $\mbf{X}^{(k)}_j,\,(j=1,\ldots,N_{k+1})$, and control, $\mbf{U}^{(k)}_j,\,(j=1,\ldots,N_k)$, in mesh interval $\mcl{S}_k,\,(k=1,\ldots,K)$.
Given this discrete approximation, the objective is to obtain an error estimate in each mesh interval based on the differences between the Lagrange polynomial approximation of the state in Eq.~\eqref{eq:lgr-xappx} and state approximations obtained via {\em{forward}} and {\em{backward}} time-marching.
In order to approximate the state via forward time-marching, an IVP must be formulated and solved.
In mesh interval $\mcl{S}_k,\,(k=1,\ldots,K)$, the single-interval IVP is formulated as
\begin{equation}\label{eq:mesh-sim-ivp}
    \deriv{\mbfh{X}^{(k)}(\zeta)}{\zeta}=\alpha\beta_k\mbf{a}(\mbfh{X}^{(k)}(\zeta),\mbf{U}^{(k)}(\zeta),\zeta),
    \qquad
    \mbfh{X}^{(k)}(-1)=\mbf{X}^{(k)}_1.
\end{equation}
Suppose that explicitly simulating Eq.~\eqref{eq:mesh-sim-ivp} from $\zeta=-1$ to $\zeta=+1$ in every mesh interval yields approximated values of the state, $\mbfh{X}^{(k)}(\mh{\zeta}^{(k)}_j),\,(j=1,\ldots,\mh{P}_k;~k=1,\ldots,K)$, at the $\mh{P}_k$ forward propagation points $\{\mh{\zeta}^{(k)}_1,\ldots,\mh{\zeta}^{(k)}_{\mh{P}_k}\}\in[-1,+1]$.
Similarly, in order to approximate the state via backward time-marching, a TVP must be formulated and solved.
In mesh interval $\mcl{S}_k,\,(k=1,\ldots,K)$, the single-interval TVP is formulated as
\begin{equation}\label{eq:mesh-sim-tvp}
    \deriv{\mbfc{X}^{(k)}(\zeta)}{\zeta}=\alpha\beta_k\mbf{a}(\mbfc{X}^{(k)}(\zeta),\mbf{U}^{(k)}(\zeta),\zeta),
    \qquad
    \mbfc{X}^{(k)}(+1)=\mbf{X}^{(k)}_{N_k+1}.
\end{equation}
Suppose that explicitly simulating Eq.~\eqref{eq:mesh-sim-tvp} from $\zeta=+1$ to $\zeta=-1$ in every mesh interval yields approximated values of the state, $\mbfc{X}^{(k)}(\mc{\zeta}^{(k)}_l),\,(l=1,\ldots,\mc{P}_k;~k=1,\ldots,K)$, at the $\mc{P}_k$ backward propagation points $\{\mc{\zeta}^{(k)}_1,\ldots,\mc{\zeta}^{(k)}_{\mc{P}_k}\}\in[-1,+1]$.
It is noted that the sets of points $\{\zeta^{(k)}_1,\ldots,\zeta^{(k)}_{N_k+1}\}$, $\{\mh{\zeta}^{(k)}_1,\ldots,\mh{\zeta}^{(k)}_{\mh{P}_k}\}$, and $\{\mc{\zeta}^{(k)}_1,\ldots,\mc{\zeta}^{(k)}_{\mc{P}_k}\}$ are not all necessarily equivalent in mesh interval $\mcl{S}_k$.
In order to obtain solutions to Eqs.~\eqref{eq:mesh-sim-ivp} and~\eqref{eq:mesh-sim-tvp}, a control function, $\mbf{U}^{(k)}(\zeta),\,(k=1,\ldots,K)$, is required because time-marching schemes take integration steps to points other than the collocation points.
In this work, a Lagrange polynomial approximation of the control is defined in mesh interval $\mcl{S}_k$ as
\begin{equation}\label{eq:mesh-sim-uappx}
    \mbf{U}^{(k)}(\zeta)=\sum_{j=1}^{N_k}\hat{\ell}_j^{(k)}(\zeta)\mbf{U}_j^{(k)},
    \qquad
    \mh{\ell}_j^{(k)}(\zeta)=\prod_{
    \begin{subarray}{c}
        l=1 \\
        l\neq{j} \\
    \end{subarray}
    }^{N_k}\dfrac{\zeta-\zeta_l^{(k)}}{\zeta_j^{(k)}-\zeta_l^{(k)}},
    \qquad
    (k=1,\ldots,K),
\end{equation}
where $\{\zeta_1^{(k)},\ldots,\zeta_{N_k}^{(k)}\}\in[-1,+1)$ is the set of $N_k$ LGR collocation points.
Equation~\eqref{eq:mesh-sim-uappx} differs from the state approximation in Eq.~\eqref{eq:lgr-xappx} as it does not include the noncollocated point.
Finally, any numerical method for solving a set of ODEs (for example, \ode{45} or \ode{89} in \MATLAB) can be employed to solve Eqs.~\eqref{eq:mesh-sim-ivp} and~\eqref{eq:mesh-sim-tvp}.
Using any such ODE solver, the sets of propagation points can either be specified by the user or determined by the ODE solver.
Performance of the mesh refinement method described in this paper is impacted by the ODE solver and accuracy tolerance of the solver.

% ------------------------------------------------------------------------------------------------ %
% ------------------- SUBSECTION: Numerical Simulation in a Mesh Interval (END) ------------------ %
% ------------------------------------------------------------------------------------------------ %

% ------------------------------------------------------------------------------------------------ %
% ---------------- SUBSECTION: Relative Error Estimate in a Mesh Interval (START) ---------------- %
% ------------------------------------------------------------------------------------------------ %
\subsection{Relative Error Estimate in a Mesh Interval}\label{subsec:mesh-err}

In this section, an estimate of the relative error in the solution on a given mesh is derived.
In the LGR collocation method, a uniquely defined function approximation is only available for the state; therefore, a relative error estimate is developed for the state only.
The key idea is that numerically solving a set of ODEs via collocation and/or time-marching schemes should yield nearly identical results, barring any propagation error or discrepancy between the accuracy of the numerical methods.
Thus, the differences between the Lagrange polynomial approximation of the state obtained by solving the NLP of Eqs.~\eqref{eq:lgr-obj}--\eqref{eq:lgr-cnt} and state approximations obtained by solving the IVP and TVP of Eqs.~\eqref{eq:mesh-sim-ivp} and \eqref{eq:mesh-sim-tvp} via forward and backward time-marching, respectively, should yield an approximation of the error in the state.
Consequently, using this newly developed error estimate ensures that the solutions obtained via collocation and explicit simulation schemes are in agreement with one another on the final mesh.

Suppose that the values of the state approximation given in Eq.~\eqref{eq:lgr-xappx} in mesh interval $\mcl{S}_k,\,(k=1,\ldots,K)$, at the forward and backward propagation points $\{\mh{\zeta}^{(k)}_1,\ldots,\mh{\zeta}^{(k)}_{\mh{P}_k}\}$ and $\{\mc{\zeta}^{(k)}_1,\ldots,\mc{\zeta}^{(k)}_{\mc{P}_k}\}$ are denoted $\{\mbf{X}^{(k)}(\mh{\zeta}^{(k)}_1),\ldots,\mbf{X}^{(k)}(\mh{\zeta}^{(k)}_{\mh{P}_k})\}$ and $\{\mbf{X}^{(k)}(\mc{\zeta}^{(k)}_1),\ldots,\mbf{X}^{(k)}(\mc{\zeta}^{(k)}_{\mc{P}_k})\}$, respectively.
Then, the forward and backward {\em{relative errors}} in the $i^{\text{th}}$ component of the state, $(i=1,\ldots,n_x)$, at the points $\{\mh{\zeta}^{(k)}_1,\ldots,\mh{\zeta}^{(k)}_{\mh{P}_k}\}$ and $\{\mc{\zeta}^{(k)}_1,\ldots,\mc{\zeta}^{(k)}_{\mc{P}_k}\}$ are defined, respectively, as
\begin{equation}\label{eq:mesh-err-rel}
    \begin{array}{rcl}
        \mh{e}^{(k)}_i(\mh{\zeta}^{(k)}_j) & = & \gamma_i\left|\mh{X}^{(k)}_i(\mh{\zeta}^{(k)}_j)-X^{(k)}_i(\mh{\zeta}^{(k)}_j)\right|, \\[0.75em]
        \mc{e}^{(k)}_i(\mc{\zeta}^{(k)}_l) & = & \gamma_i\left|\mc{X}^{(k)}_i(\mc{\zeta}^{(k)}_l)-X^{(k)}_i(\mc{\zeta}^{(k)}_l)\right|, \\
    \end{array}
    \qquad
    \left(
    \begin{array}{ll}
        j & =1,\ldots,\mh{P}_k \\
        k & =1,\ldots,K \\
        l & =1,\ldots,\mc{P}_k \\
    \end{array}
    \right)\!,
\end{equation}
where the corresponding {\em{error scaling factor}}, $\gamma_i\in\mbb{R}$, is defined such that
\begin{equation}\label{eq:mesh-err-gamma}
    \gamma_i^{-1}=1+\underset{k\in\{1,\ldots,K\}}{\underset{j\in\{1,\ldots,N_k+1\}}{\max}}\left|X^{(k)}_i(\zeta^{(k)}_j)\right|,
    \qquad
    (i=1,\ldots,n_x),
\end{equation}
and $X^{(k)}_i(x)$ denotes the $i^{\text{th}}$ component of the state approximation $\mbf{X}^{(k)}(\zeta)$ evaluated at the specified point (in this case, $x$).
To account for both relative errors in Eq.~\eqref{eq:mesh-err-rel}, the {\em{maximum relative error}} in mesh interval $\mcl{S}_k,\,(k=1,\ldots,K)$, is defined as
\begin{equation}\label{eq:mesh-err-relmax}
    e^{(k)}_{\max}=\underset{i\in\{1,\ldots,n_x\}}{\max}\lprp{\,
    \underset{j\in\{1,\ldots,\mh{P}_k\}}{\max}\mh{e}^{(k)}_i(\mh{\zeta}^{(k)}_j),
    \underset{l\in\{1,\ldots,\mc{P}_k\}}{\max}\mc{e}^{(k)}_i(\mc{\zeta}^{(k)}_l)}\!.
\end{equation}
% ------------------------------------------------------------------------------------------------ %
% ----------------- SUBSECTION: Relative Error Estimate in a Mesh Interval (END) ----------------- %
% ------------------------------------------------------------------------------------------------ %

% ------------------------------------------------------------------------------------------------ %
% -------------------- SUBSECTION: Mesh Refinement in a Mesh Interval (START) -------------------- %
% ------------------------------------------------------------------------------------------------ %
\subsection{Mesh Refinement in a Mesh Interval}\label{subsec:mesh-refine}

The objective of the mesh refinement method is to meet the desired mesh tolerance, $\epsilon$, in every mesh interval on the smallest mesh (that is, with the fewest total number of collocation points).
For simplicity, let $\mcl{K}^{+},\,(\mcl{K}^{+}\subseteq\{1,\ldots,K\})$, denote the set of mesh intervals in which the desired mesh tolerance is not met (that is, $e^{(k)}_{\max}>\epsilon,\,(k\in\mcl{K}^{+})$); furthermore, let $\mcl{K}^{-},\,(\mcl{K}^{-}\subseteq\{1,\ldots,K\})$, denote the set of mesh intervals in which the desired mesh tolerance is met (that is, $e^{(k)}_{\max}\leq\epsilon,\,(k\in\mcl{K}^{-})$).
It is noted that $\mcl{K}^{+}\cap\mcl{K}^{-}=\emptyset$ and $\mcl{K}^{+}\cup\mcl{K}^{-}=\{1,\ldots,K\}$.
The current mesh is refined only if $\mcl{K}^{+}\neq\emptyset$.
First, the method attempts to increase the degree of the polynomial approximation in a mesh interval and/or divide a mesh interval into subintervals.
Then, the method attempts to merge adjacent mesh intervals and/or decrease the degree of the polynomial approximation in a mesh interval.
For use hereafter, let the notation $x^{[M]}$ denote a quantity (in this case, $x$) on mesh refinement iteration $M$, where $M=0$ corresponds to the initial mesh.

% ----- SUBSUBSECTION: Method for Increasing the Polynomial Approximation Degree ... (START) ----- %
\subsubsection{\ensuremath{p} Refinement: Increasing the Polynomial Degree in a Mesh Interval}\label{subsubsec:mesh-refine-pp}

To reduce the maximum relative error, $p$ refinement attempts to strictly increase the number of collocation points $N^{[M]}_k$ in mesh interval $\mcl{S}_k,\,(k\in\mcl{K}^{+})$, on mesh $M$ to
\begin{equation}\label{eq:mesh-refine-pp-Nk}
    N_{k}^{[M+1]}=N_{k}^{[M]}+P^{+}_k,
    \qquad
    P^{+}_k=\lcrc{\log_{10}\left(\dfrac{e_{\max}^{(k)}}{\epsilon}\right)}\!,
\end{equation}
where $N_k^{[M+1]}$ denotes the number of collocation points in $\mcl{S}_k$ on the ensuing mesh, $\lcrc{\cdot}$ denotes the ceiling function, and $P^{+}_k\geq{1}$ because $e_{\max}^{(k)}>\epsilon$.
Equation~\eqref{eq:mesh-refine-pp-Nk} removes the arbitrary constant found in the $p$ refinement step of Ref.~\cite{DarbyRao2011a} in order to keep the size of the mesh as small as possible.
To ensure that the approximating polynomial degree in a mesh interval does not grow to an unreasonably large value, a user-specified upper limit $N_{\max}\geq{2}$ is set for the maximum allowable polynomial approximation degree.
If $N_{\max}$ is exceeded (that is, $p$ refinement is exhausted), then the mesh interval is divided into subintervals using the method presented in Section~\ref{subsubsec:mesh-refine-hp}.

% ------ SUBSUBSECTION: Method for Increasing the Polynomial Approximation Degree ... (END) ------ %

% ----------- SUBSUBSECTION: Method for Dividing a Mesh Interval (h Refinement) (START) ---------- %
\subsubsection{\ensuremath{h} Refinement: Dividing a Mesh Interval into Subintervals}\label{subsubsec:mesh-refine-hp}

Suppose that evaluating Eq.~\eqref{eq:mesh-refine-pp-Nk} yields $N_{k}^{[M+1]}>N_{\max}$, which indicates that mesh interval $\mcl{S}_k,\,(k\in\mcl{K}^{+})$, on mesh $M$ must be divided into subintervals.
Identical to the $h$ refinement method described in Ref.~\cite{PattersonRao2015}, it is desired to keep the predicted number of total collocation points on the ensuing mesh and employ $N_{\min}$ collocation points in each newly created subinterval, where $N_{\min}$ denotes a user-specified minimum allowable polynomial approximation degree such that $N_{\max}\geq{N_{\min}}\geq{2}$.
Thus, mesh interval $\mcl{S}_k,\,(k\in\mcl{K}^{+})$, is divided into $H_k$ uniformly spaced subintervals on the ensuing mesh, where
\begin{equation}\label{eq:mesh-refine-hp-Hk}
    H_k=\max\lprp{2,\lcrc{\dfrac{N_{k}^{[M+1]}}{N_{\min}}}}\!.
\end{equation}
The new subintervals $\mh{\mcl{S}}_q=[\mh{\tau}^{(k)}_{q-1},\mh{\tau}^{(k)}_q]\subset\mcl{S}_k=[\tau_{k-1},\tau_k],\,(q=1,\ldots,H_k;~k\in\mcl{K}^{+})$, are created by defining the mesh points on the ensuing mesh as
\begin{equation}\label{eq:mesh-refine-hp-tauq}
    \mh{\tau}^{(k)}_q=\tau_{k-1}+q\lprp{\dfrac{\tau_k-\tau_{k-1}}{H_k}}\!,
    \qquad
    (q=0,\ldots,H_k),
\end{equation}
where $\tau_{k-1}=\mh{\tau}^{(k)}_0<\mh{\tau}^{(k)}_1<\cdots<\mh{\tau}^{(k)}_{H_k-1}<\mh{\tau}^{(k)}_{H_k}=\tau_k$.

% ------------ SUBSUBSECTION: Method for Dividing a Mesh Interval (h Refinement) (END) ----------- %

% -------- SUBSUBSECTION: Method for Merging Adjacent Mesh Intervals (h Reduction) (START) ------- %
\subsubsection{\ensuremath{h} Reduction: Merging Adjacent Mesh Intervals}\label{subsubsec:mesh-refine-hm}

In addition to increasing the mesh size via $p$ and $h$ refinement, the mesh size can be decreased via either $h$ or $p$ reduction.
The method developed in this research for $h$ reduction utilizes state approximations obtained via forward and backward time marching over the domain of both mesh intervals to determine if merging is possible, which extends the single-interval ideas of Sections~\ref{subsec:mesh-sim} and~\ref{subsec:mesh-err} to adjacent mesh intervals.
It is noted that this $h$ reduction step differs significantly from that described in Ref.~\cite{LiuRao2015}, where different polynomial approximations of the state are obtained from power series expansions about the shared mesh point.

First, let $\mcl{Q},\,(\mcl{Q}\subseteq\mcl{K}^{-};~K\notin\mcl{Q})$, denote the set of adjacent mesh interval pairs that both satisfy the desired mesh tolerance.
The goal is to merge adjacent mesh intervals $\mcl{S}_k$ and $\mcl{S}_{k+1},\,(\{k,k+1\}\in\mcl{K}^{-})$, to form a single mesh interval $\mb{\mcl{S}}_q=\mcl{S}_{k}\cup\mcl{S}_{k+1}=[\tau_{k-1},\tau_{k+1}],\,(q\in\mcl{Q})$.
For simplicity, let the following adjacent mesh interval transformations be
\begin{equation}\label{eq:mesh-refine-hm-xichi}
    \begin{array}{rclcl}
        \chi(x) & \equiv & \chi(x,\tau_{k-1},\tau_k,\tau_{k+1}) & = & \zeta(\tau(x,\tau_{k-1},\tau_k),\tau_k,\tau_{k+1}), \\
        \xi(x)  & \equiv & \xi(x,\tau_{k-1},\tau_k,\tau_{k+1})  & = & \zeta(\tau(x,\tau_k,\tau_{k+1}),\tau_{k-1},\tau_k), \\
    \end{array}
    \qquad
    (k=1,\ldots,K-1),
\end{equation}
where $\chi(x)$ transforms a quantity (in this case, $x$) relative to $\zeta^{(k)}\in[-1,+1]$ to its corresponding value relative to $\zeta^{(k+1)}\in[-1,+1]$, and $\xi(x)$ provides the inverse transformation.
It is noted that Eq.~\eqref{eq:mesh-refine-hm-xichi} utilizes the transformations in Eq.~\eqref{eq:bolza-zeta2tau_tau2zeta}.
Similar to Eqs.~\eqref{eq:mesh-sim-ivp} and \eqref{eq:mesh-sim-tvp}, an IVP and TVP must both be formulated and solved in order to approximate the state over the domain of the adjacent mesh intervals via forward and backward time-marching, which are defined as follows.
First, the multiple-interval IVP in mesh interval $\mb{\mcl{S}}_q,\,(q\in\mcl{Q})$, is formulated as
\begin{equation}\label{eq:mesh-refine-hm-ivp}
    \deriv{\mbfbh{X}^{(q)}(\zeta)}{\zeta}=\alpha\beta_k\mbf{a}(\mbfbh{X}^{(q)}(\zeta),\mbf{U}^{(k)}(\zeta),\zeta),
    \qquad
    \mbfbh{X}^{(q)}(-1)=\mbf{X}^{(k)}_1.
\end{equation}
Suppose that explicitly simulating Eq.~\eqref{eq:mesh-refine-hm-ivp} from $\zeta=-1$ to $\zeta=\xi(\zeta^{(k+1)}_{N_{k+1}+1})>+1$ (that is, forward from mesh interval $\mcl{S}_k$ into $\mcl{S}_{k+1}$) yields approximated values of the state, $\mbfbh{X}^{(q)}(\mbh{\zeta}^{(q)}_j),\,(j=1,\ldots,\mbh{P}_q)$, at the $\mbh{P}_q$ extended forward propagation points $\{\mbh{\zeta}^{(q)}_1,\ldots,\mbh{\zeta}^{(q)}_{\mbh{P}_q}\}\in[-1,\xi(\zeta^{(k+1)}_{N_{k+1}+1})]$ in mesh interval $\mb{\mcl{S}}_q,\,(q\in\mcl{Q})$.
Next, the multiple-interval TVP in mesh interval $\mb{\mcl{S}}_q,\,(q\in\mcl{Q})$, is formulated as
\begin{equation}\label{eq:mesh-refine-hm-tvp}
    \deriv{\mbfbc{X}^{(q)}(\zeta)}{\zeta}=\alpha\beta_{k+1}\mbf{a}(\mbfbc{X}^{(q)}(\zeta),\mbf{U}^{(k+1)}(\zeta),\zeta),
    \qquad
    \mbfbc{X}^{(q)}(+1)=\mbf{X}^{(k+1)}_{N_{k+1}+1}.
\end{equation}
Suppose that explicitly simulating Eq.~\eqref{eq:mesh-refine-hm-tvp} from $\zeta=+1$ to $\zeta=\chi(\zeta^{(k)}_1)<-1$ (that is, backward from mesh interval $\mcl{S}_{k+1}$ into $\mcl{S}_k$) yields approximated values of the state, $\mbfbc{X}^{(q)}(\mbc{\zeta}^{(q)}_l),\,(l=1,\ldots,\mbc{P}_q)$, at the $\mbc{P}_q$ extended backward propagation points $\{\mbc{\zeta}^{(q)}_1,\ldots,\mbc{\zeta}^{(q)}_{\mbc{P}_q}\}\in[\chi(\zeta^{(k)}_1),+1]$ in mesh interval $\mb{\mcl{S}}_q,\,(q\in\mcl{Q})$.
It is noted that the control in mesh interval $\mcl{S}_k$ is extrapolated forward into mesh interval $\mcl{S}_{k+1}$ in the multiple-interval IVP of Eq.~\eqref{eq:mesh-refine-hm-ivp}, whereas the control in mesh interval $\mcl{S}_{k+1}$ is extrapolated backward into mesh interval $\mcl{S}_k$ in the multiple-interval TVP of Eq.~\eqref{eq:mesh-refine-hm-tvp}.

In order to determine whether or not the adjacent mesh intervals can be merged, the relative error is estimated over the domain of the adjacent mesh intervals using the state approximations obtained by solving the multiple-interval IVP and TVP of Eqs.~\eqref{eq:mesh-refine-hm-ivp} and \eqref{eq:mesh-refine-hm-tvp}, respectively.
Because the state approximations given by Eq.~\eqref{eq:lgr-xappx} in the adjacent mesh intervals are not necessarily equivalent, the appropriate state approximations must be compared at the corresponding propagation points.
Similar to Eq.~\eqref{eq:mesh-err-rel}, the relative errors in the $i^{\text{th}}$ component of the state at the extended forward and backward propagation points $\{\mbh{\zeta}^{(q)}_1,\ldots,\mbh{\zeta}^{(q)}_{\mbh{P}_q}\}$ and $\{\mbc{\zeta}^{(q)}_1,\ldots,\mbc{\zeta}^{(q)}_{\mbc{P}_q}\}$ in mesh interval $\mb{\mcl{S}}_q,\,(q\in\mcl{Q})$, are defined, respectively, as
\begin{equation}\label{eq:mesh-refine-hm-rel}
    \renewcommand*{\mainterm}{\gamma_i\left|\mbh{X}^{(q)}_i(\mbh{\zeta}^{(q)}_j)-X^{(k+1)}_i(\chi(\mbh{\zeta}^{(q)}_j))\right|}
    \renewcommand*{\leftterm}{\chi(\zeta^{(k)}_1)}
    \renewcommand*{\centerterm}{\mbh{\zeta}^{(q)}_j}
    \renewcommand*{\rightterm}{\xi(\zeta^{(k+1)}_{N_{k+1}+1})}
    \settowidth{\mainside}{\term{\mainterm}~}
    \settowidth{\leftside}{\term{\leftterm}}
    \settowidth{\centerside}{\term{\centerterm}}
    \settowidth{\rightside}{\term{\rightterm}}
    \begin{array}{rcl}
        \mbh{e}^{(q)}_i(\mbh{\zeta}_j^{(q)}) & = & \left\{
        \begin{aligned}
            \makebox[\mainside][l]{\term{\gamma_i\left|\mbh{X}^{(q)}_i(\mbh{\zeta}^{(q)}_j)-X^{(k)}_i(\mbh{\zeta}^{(q)}_j)\right|}}
            & ~\text{for} & \makebox[\leftside][r]{\term{-1}} & \leq & \makebox[\centerside][c]{\term{\mbh{\zeta}^{(q)}_j}} & \leq & \makebox[\rightside][l]{\term{+1,}} \\[0.5em]
            \makebox[\mainside][l]{\term{\gamma_i\left|\mbh{X}^{(q)}_i(\mbh{\zeta}^{(q)}_j)-X^{(k+1)}_i(\chi(\mbh{\zeta}^{(q)}_j))\right|}}
            & ~\text{for} & \makebox[\leftside][r]{\term{+1}} & <    & \makebox[\centerside][c]{\term{\mbh{\zeta}^{(q)}_j}} & \leq & \makebox[\rightside][l]{\term{\xi(\zeta^{(k+1)}_{N_{k+1}+1}),}} \\
        \end{aligned}
        \right. \\[2.0em]
        \mbc{e}^{(q)}_i(\mbc{\zeta}_l^{(q)}) & = &
        \left\{
        \begin{aligned}
            \makebox[\mainside][l]{\term{\gamma_i\left|\mbc{X}^{(q)}_i(\mbc{\zeta}^{(q)}_l)-X^{(k)}_i(\xi(\mbc{\zeta}^{(q)}_l))\right|}}
            & ~\text{for} & \makebox[\leftside][r]{\term{\chi(\zeta^{(k)}_1)}} & \leq & \makebox[\centerside][c]{\term{\mbc{\zeta}^{(q)}_l}} & <    & \makebox[\rightside][l]{\term{-1,}} \\[0.5em]
            \makebox[\mainside][l]{\term{\gamma_i\left|\mbc{X}^{(q)}_i(\mbc{\zeta}^{(q)}_l)-X^{(k+1)}_i(\mbc{\zeta}^{(q)}_l)\right|}}
            & ~\text{for} & \makebox[\leftside][r]{\term{-1}}                  & \leq & \makebox[\centerside][c]{\term{\mbc{\zeta}^{(q)}_l}} & \leq & \makebox[\rightside][l]{\term{+1,}} \\
        \end{aligned}
        \right. \\[2.0em]
        \multicolumn{3}{c}{(i=1,\ldots,n_x;~j=1,\ldots,\mbh{P}_q;~\{k,k+1\}\in\mcl{K}^{-};~l=1,\ldots,\mbc{P}_q),} \\
    \end{array}
\end{equation}
where $X^{(k)}_i(x)$ and $X^{(k+1)}_i(x)$ denote the $i^{\text{th}}$ component of the respective state approximations $\mbf{X}^{(k)}(\zeta)$ and $\mbf{X}^{(k+1)}(\zeta)$ given by Eq.~\eqref{eq:lgr-xappx} evaluated at the specified point (in this case, $x$).
To account for both relative errors in Eq.~\eqref{eq:mesh-refine-hm-rel}, the maximum relative error in mesh interval $\mb{\mcl{S}}_q,\,(q\in\mcl{Q})$, is defined as
\begin{equation}\label{eq:mesh-refine-hm-relmax}
    \mb{e}^{(q)}_{\max}=\underset{i\in\{1,\ldots,n_x\}}{\max}\lprp{\,
    \underset{j\in\{1,\ldots,\mbh{P}_q\}}{\max}\mbh{e}^{(q)}_i(\mbh{\zeta}^{(q)}_j),
    \underset{l\in\{1,\ldots,\mbc{P}_q\}}{\max}\mbc{e}^{(q)}_i(\mbc{\zeta}^{(q)}_l)}\!.
\end{equation}
It is noted that the adjacent mesh intervals cannot be merged if $\mb{e}^{(q)}_{\max}>\epsilon,\,(q\in\mcl{Q})$.

Once the method determines all pairs of adjacent mesh intervals that can be merged, merging then occurs in ascending order based on the overall maximum relative error obtained via Eqs.~\eqref{eq:mesh-err-relmax} and~\eqref{eq:mesh-refine-hm-relmax} (that is, based on $\max(e^{(k)}_{\max},e^{(k+1)}_{\max},\mb{e}^{(q)}_{\max})$).
If $h$ reduction is performed, the number of collocation points $N^{[M+1]}_q$ in mesh interval $\mb{\mcl{S}}_q,\,(q\in\mcl{Q})$, on mesh $M+1$ is determined by
\begin{equation}\label{eq:mesh-refine-hm-Nk}
    N^{[M+1]}_q=\max\lprp{N^{[M]}_k,N^{[M]}_{k+1}},
\end{equation}
where the mesh point $\tau_k$ is removed on the ensuing mesh.
It is important to note that the $h$ reduction step presented in this research does not require the number of collocation points in adjacent mesh intervals to be equal, which can be a limiting factor enforced by the method described in Ref.~\cite{LiuRao2015}.

% --------- SUBSUBSECTION: Method for Merging Adjacent Mesh Intervals (h Reduction) (END) -------- %

% ----- SUBSUBSECTION: Method for Decreasing the Polynomial Approximation Degree ... (START) ----- %
\subsubsection{\ensuremath{p} Reduction: Decreasing the Polynomial Degree in a Mesh Interval}\label{subsubsec:mesh-refine-pm}

Before attempting $p$ reduction, all possible mesh intervals are merged using the method described in Section~\ref{subsubsec:mesh-refine-hm}.
If mesh interval $\mcl{S}_k,\,(k\in\mcl{K}^{-})$, is not merged, then $p$ reduction attempts to decrease the number of collocation points $N^{[M]}_k$ on mesh $M$ to
\begin{equation}\label{eq:mesh-refine-pm-Nk}
    \begin{array}{cc}
        N_{k}^{[M+1]}=\max\lprp{N_{\min},N_{k}^{[M]}-P^{-}_k},
        \qquad
        P^{-}_k=\lfrf{\log_{10}\lprp{\dfrac{\epsilon}{e_{\max}^{(k)}}}^{1/\delta}}\!, \\
        \delta=N_{\min}+N_{\max}-N^{[M]}_k,
    \end{array}
\end{equation}
where $\lfrf{\cdot}$ denotes the floor function, and $P^{-}_k\geq{0}$ because $e_{\max}^{(k)}\leq\epsilon$.
The parameter $\delta$ in Eq.~\eqref{eq:mesh-refine-pm-Nk} controls the acceptable order of magnitude difference between the current maximum relative error in a mesh interval and the desired mesh tolerance, where a larger number of collocation points already present in a mesh interval yields a smaller value for $\delta$ and vice versa.
When $P^{-}_k=0$, the polynomial approximation degree in mesh interval $\mcl{S}_k,\,(k\in\mcl{K}^{-})$, cannot be reduced.
Different from the $p$ reduction step described in Ref.~\cite{LiuRao2015} that utilizes a different polynomial approximation of the state obtained from a power series expansion about the mesh interval midpoint, the $p$ reduction step of Eq.~\eqref{eq:mesh-refine-pm-Nk} follows similar logic to that employed in Eq.~\eqref{eq:mesh-refine-pp-Nk}.
Equation~\eqref{eq:mesh-refine-pm-Nk} also ensures that the number of collocation points in a mesh interval is not decreased below the lower limit $N_{\min}$.

% ------ SUBSUBSECTION: Method for Decreasing the Polynomial Approximation Degree ... (END) ------ %

% ------------------------------------------------------------------------------------------------ %
% --------------------- SUBSECTION: Mesh Refinement in a Mesh Interval (END) --------------------- %
% ------------------------------------------------------------------------------------------------ %

% ------------------------------------------------------------------------------------------------ %
% ------------------ SUBSECTION: Adaptive Mesh Refinement Method Summary (START) ----------------- %
% ------------------------------------------------------------------------------------------------ %
\subsection{Summary of Adaptive Mesh Refinement Method}\label{subsec:mesh-sum}

A summary of the method developed in this research is shown in Algorithm~\ref{alg:mesh-sum}, which consists of three major components: numerical simulation, relative state error estimation, and mesh refinement in a mesh interval, as discussed in Sections~\ref{subsec:mesh-sim}--\ref{subsec:mesh-refine}, respectively.
% ------------------------------------------
% Algorithm: adaptive mesh refinement method
\begin{algorithm}[!t]
    \scriptsize
    \caption{\enskip{Adaptive Mesh Refinement Method}}
    \label{alg:mesh-sum}
    \begin{algorithmic}[1]
        \State Supply an initial mesh that satisfies Eq.~\eqref{eq:bolza-mesh}.
            \Comment{Section~\ref{sec:bolza}}
        \State Assign user-specified parameters: \texttt{s}, $\epsilon_{\text{ODE}}$, $\epsilon$, $M_{\max}$, $N_{\min}$, $N_{\max}$.
            \Comment{Section~\ref{subsec:mesh-sum}}
        \State $M\gets{0}$.
        \While{$M<M_{\max}$}
            \State $K\gets$ number of mesh intervals on mesh $M$.
            \State $N_k\gets$ number of collocation points in mesh interval $\mcl{S}_k$ on mesh $M$.
            \State Solve the NLP on mesh $M$ as follows:
                \Comment{Section~\ref{sec:lgr}}
            \State $\mbf{X}_{1:N_k+1}^{(k)},\mbf{U}_{1:N_k}^{(k)}~\forall~k=1\,\ldots,K\gets$ Solve Eqs.~\eqref{eq:lgr-obj}--\eqref{eq:lgr-cnt}.
            \For{$k=1,\ldots,K$}
                \State Explicitly simulate the dynamics as follows:
                    \Comment{Section~\ref{subsec:mesh-sim}}
                \State $\mbfh{X}_{1:\mh{P}_k}^{(k)}\gets$ Solve the IVP of Eq.~\eqref{eq:mesh-sim-ivp} in $\mcl{S}_k$.
                \State $\mbfc{X}_{1:\mc{P}_k}^{(k)}\gets$ Solve the TVP of Eq.~\eqref{eq:mesh-sim-tvp} in $\mcl{S}_k$.
                \State Compute the maximum relative error estimate as follows:
                    \Comment{Section~\ref{subsec:mesh-err}}
                \State $e^{(k)}_{\max}\gets$ Solve Eq.~\eqref{eq:mesh-err-relmax} in $\mcl{S}_k$.
            \EndFor
            \If{$e^{(k)}_{\max}\leq\epsilon~\forall~k=1,\ldots,K$}
                \State \textbf{exit}
            \EndIf
            \State Refine mesh $M$ as follows:
                \Comment{Section~\ref{subsec:mesh-refine}}
            \For{$k=1,\ldots,K$}
                \If{$e^{(k)}_{\max}>\epsilon$}
                    \State Attempt $p$ refinement in $\mcl{S}_k$.
                        \Comment{Section~\ref{subsubsec:mesh-refine-pp}}
                    \If{$p$ refinement fails in $\mcl{S}_k$}
                        \State Perform $h$ refinement in $\mcl{S}_k$.
                            \Comment{Section~\ref{subsubsec:mesh-refine-hp}}
                    \EndIf
                \Else
                    \If{$k<K$ \textbf{and} $e^{(k+1)}_{\max}\leq\epsilon$}
                        \State Attempt $h$ reduction in $\mcl{S}_k$ and $\mcl{S}_{k+1}$.
                            \Comment{Section~\ref{subsubsec:mesh-refine-hm}}
                        \If{$h$ reduction fails in $\mcl{S}_k$ and $\mcl{S}_{k+1}$}
                            \State Attempt $p$ reduction in $\mcl{S}_k$.
                                \Comment{Section~\ref{subsubsec:mesh-refine-pm}}
                        \EndIf
                    \Else
                        \State Attempt $p$ reduction in $\mcl{S}_k$.
                            \Comment{Section~\ref{subsubsec:mesh-refine-pm}}
                    \EndIf
                \EndIf
            \EndFor
            \State $M\gets{M+1}$.
        \EndWhile
    \end{algorithmic}
    \normalsize
\end{algorithm}%
% ------------------------------------------
Mesh refinement in a mesh interval consists of four steps: $p$ refinement, $h$ refinement, $h$ reduction, and $p$ reduction, as discussed in Sections~\ref{subsubsec:mesh-refine-pp}--\ref{subsubsec:mesh-refine-pm}, respectively.
Six user-specified parameters are required: an ODE solver, $\texttt{s}$, an ODE solver tolerance, $\epsilon_{\text{ODE}}$, a mesh tolerance, $\epsilon$, a maximum number of mesh refinement iterations, $M_{\max}$, and minimum and maximum allowable numbers of collocation points, $N_{\min}$ and $N_{\max}$, respectively.
As with any mesh refinement method, performance of the method depends on the method parameters, initial mesh, and mesh accuracy tolerance.

% ------------------------------------------------------------------------------------------------ %
% ------------------- SUBSECTION: Adaptive Mesh Refinement Method Summary (END) ------------------ %
% ------------------------------------------------------------------------------------------------ %

% ------------------------------------------------------------------------------------------------ %
% ------------------------------------------------------------------------------------------------ %
% ------------------------ SECTION: Adaptive Mesh Refinement Method (End) ------------------------ %
% ------------------------------------------------------------------------------------------------ %
% ------------------------------------------------------------------------------------------------ %

% ------------------------------------------------------------------------------------------------ %
% ------------------------------------------------------------------------------------------------ %
% ----------------------------------- SECTION: Examples (START) ---------------------------------- %
% ------------------------------------------------------------------------------------------------ %
% ------------------------------------------------------------------------------------------------ %
\section{Examples}\label{sec:examples}

In this section, the mesh refinement method of Section~\ref{sec:mesh} is demonstrated on three examples from the open literature.
The first example is the minimum-time supersonic aircraft climb problem taken from Ref.~\cite{DarbyRao2011a} that uses the dynamic model provided in Ref.~\cite{SeywaldWell1994}.
This first example demonstrates how the method is able to handle a challenging real-world problem with an active state constraint and a smooth control profile.
The second example is the minimum-time robot arm reorientation problem taken from Ref.~\cite{DolanMunson2004}.
This second example demonstrates how the method is able to handle multiple discontinuities in the control.
The third example is the hyper-sensitive optimal control problem taken from Refs.~\cite{RaoMease2000,DarbyRao2011a}.
This third example addresses how the method may need to be modified when explicit simulation fails in either forward or backward simulation.
All three examples demonstrate, to various extents, the ability of the method to effectively converge to an optimal solution on a smaller mesh that satisfies the desired mesh tolerance compared to the methods of Refs.~\cite{PattersonRao2015,LiuRao2015,LiuRao2018}.

The mesh refinement method developed in this paper is hereafter referred to as the \HR~method, where \texttt{s} refers to the user-specified numerical simulation scheme.
When comparing against various methods, the terminology \PR, \LR, and \LRL~is adopted to refer to the methods of Refs.~\cite{PattersonRao2015}, \cite{LiuRao2015}, and~\cite{LiuRao2018}, respectively, implemented with the error estimate of Ref.~\cite{PattersonRao2015}.
Moreover, the terminology \PRHR, \LRHR, and \LRLHR~is adopted to refer to the methods of Refs.~\cite{PattersonRao2015}, \cite{LiuRao2015}, and~\cite{LiuRao2018}, respectively, implemented with the error estimate of Section~\ref{subsec:mesh-err}.
The previously mentioned methods all require additional parameters.
For all methods, combinations of parameter values of $N_{\min}\in\{2,3,4\}$ and $N_{\max}\in\{6,8,\ldots,14\}$ are examined.
Similar to Ref.~\cite{LiuRao2015}, the \LR~and \LRHR~methods use a parameter value of $\bar{R}=1.2$, where $\bar{R}$ is the threshold of significance for the second derivative ratio.
Similar to Ref.~\cite{LiuRao2018}, the \LRL~and \LRLHR~methods use a parameter value of $\bar{\sigma}=0.5$, where $\bar{\sigma}$ is the threshold of significance for the decay rate of the Legendre coefficients.
For the \HR, \PRHR, \LRHR, and \LRLHR~methods, the medium and high accuracy \MATLAB~ODE solvers $\texttt{s}\in\{\ode{45},\ode{89}\}$ are used to solve the single- and multiple-interval IVPs and TVPs of Sections~\ref{subsec:mesh-sim} and~\ref{subsubsec:mesh-refine-hm} and demonstrate how the performance of the \HR~method is impacted.
The \MATLAB~ODE solvers are employed using an accuracy tolerance of $\epsilon_{\text{ODE}}=1\times{10}^{-6}$ (that is, the relative and absolute error accuracy tolerances when using an ODE solver are set to $\epsilon_{\text{ODE}}$, where the error is controlled relative to the norm of the solution).
Table~\ref{tab:examples-methods} summarizes the mesh refinement methods considered in this work with their respective error estimates and variant parameters.
% ----------------------------------------------
% Table: summary of mesh refinement methods used
\begin{table}[!b]
    \centering
    \caption{Summary of the mesh refinement methods considered in this work.}
    \begin{tabular}{lcc}
        \toprule\toprule
        \multicolumn{1}{c}{Method Name} & Error Estimation Method & Mesh Refinement Algorithm \\
        \midrule\midrule
        \PRfull    & Ref.~\cite{PattersonRao2015}  & Ref.~\cite{PattersonRao2015} \\
        \PRHRfull  & Section~\ref{subsec:mesh-err} & Ref.~\cite{PattersonRao2015} \\
        \LRfull    & Ref.~\cite{PattersonRao2015}  & Ref.~\cite{LiuRao2015} \\
        \LRHRfull  & Section~\ref{subsec:mesh-err} & Ref.~\cite{LiuRao2015} \\
        \LRLfull   & Ref.~\cite{PattersonRao2015}  & Ref.~\cite{LiuRao2018} \\
        \LRLHRfull & Section~\ref{subsec:mesh-err} & Ref.~\cite{LiuRao2018} \\
        \HRfull    & Section~\ref{subsec:mesh-err} & Section~\ref{subsec:mesh-refine} \\
        \bottomrule\bottomrule
    \end{tabular}
    \label{tab:examples-methods}
\end{table}%
% ----------------------------------------------
Note that a method variant includes all necessary parameters (for example, \HRvar{2}{6}{\ode{45}} refers to a \HR~method variant with $N_{\min}=2$, $N_{\max}=6$, and $\texttt{s}=\ode{45}$).
For mesh size analysis, let the quantities $N$ and $K$ denote the total number of collocation points and mesh intervals, respectively.
For error analysis, let the quantities $e_{\max}^{[M]}$ and $\tilde{e}_{\max}^{[M]}$ denote the maximum relative errors obtained on mesh $M$ via the error estimation methods of Section~\ref{subsec:mesh-err} and Ref.~\cite{PattersonRao2015}, respectively.

All results are obtained with the \MATLAB~general-purpose optimal control software $\mbb{GPOPS}$-$\mbb{II}$~\cite{PattersonRao2014} using the NLP solver SNOPT~\cite{GillSaunders2005}.
SNOPT is employed with a tolerance of $\epsilon_{\text{NLP}}=1\times{10}^{-8}$, where first derivatives are supplied via a sparse central difference scheme.
For all methods, a mesh refinement tolerance of $\epsilon=1\times{10}^{-6}$ and a maximum number of mesh refinement iterations of $M_{\max}=40$ are used.
For all examples, the initial mesh consists of 10 uniformly spaced intervals with $N_{\min}$ collocation points per mesh interval, and the initial guess is a straight line for variables with boundary conditions at both endpoints and constant for variables with boundary conditions at only one endpoint.
Finally, all computations are performed on an Apple M3 Pro MacBook Pro running macOS Sonoma version 14.6.1 with 36 GB LPDDR5 of unified memory using \MATLAB~version R2024a (build 24.1.0.2603908).
All reported computational run times are 10-run averages of the total time spent within SNOPT and total execution time.

% ------------------------------------------------------------------------------------------------ %
% --------- SUBSECTION: Example 1: Minimum-Time to Climb of a Supersonic Aircraft (START) -------- %
% ------------------------------------------------------------------------------------------------ %
\subsection{Example 1: Minimum-Time to Climb of a Supersonic Aircraft}\label{subsec:examples-supersonic}

\begin{figure}[!t]
    \centering
    \hspace*{\fill}
    \subfloat[$h(t)$ vs. $\tau$.]{%
        \includegraphics[height=\figureheight]{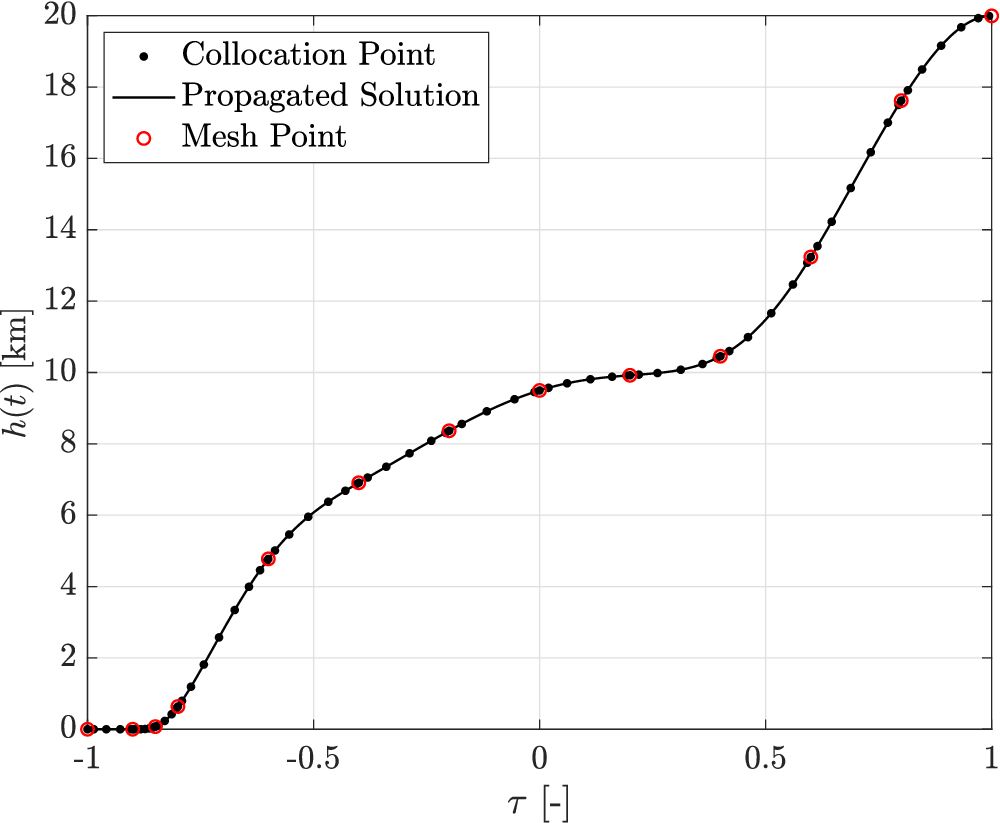}%
        \label{fig:examples-supersonic-h}
        }
    \hspace*{\fill}
    \subfloat[$v(t)$ vs. $\tau$.]{%
        \includegraphics[height=\figureheight]{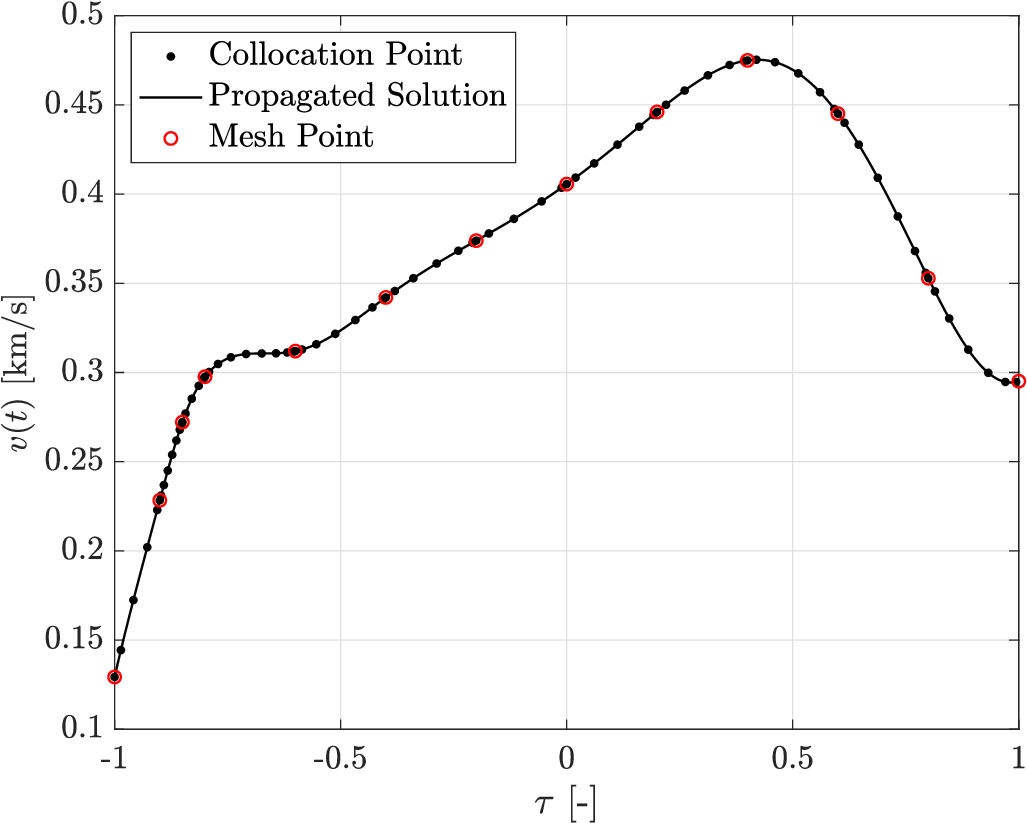}%
        \label{fig:examples-supersonic-v}
        }
    \hspace*{\fill}
    \\
    \hspace*{\fill}
    \subfloat[$\gamma(t)$ vs. $\tau$.]{%
        \includegraphics[height=\figureheight]{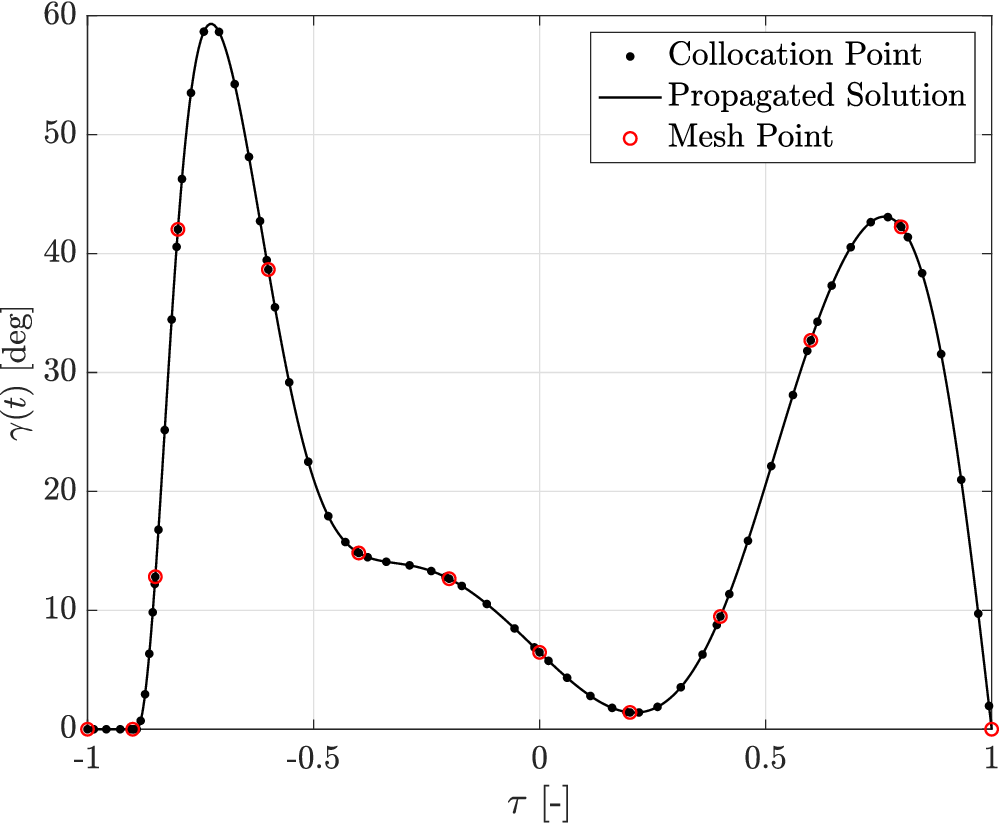}%
        \label{fig:examples-supersonic-gamma}
        }
    \hspace*{\fill}
    \caption{Optimal state components for Example 1 using the \HRvar{3}{10}{\ode{45}} method.}
    \label{fig:examples-supersonic-state}
\end{figure}

\begin{figure}[!t]
    \centering
    \hspace*{\fill}
    \includegraphics[height=\figureheight]{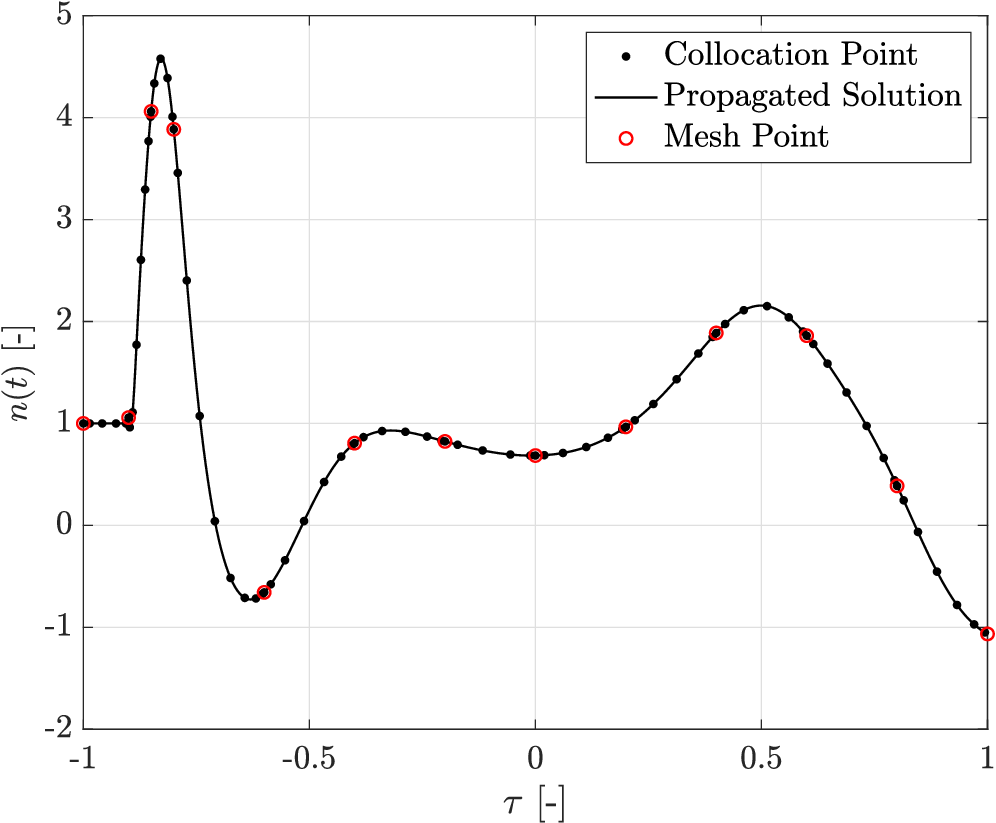}%
    \hspace*{\fill}
    \caption{Optimal control, $n(t)$ vs. $\tau$, for Example 1 using the \HRvar{3}{10}{\ode{45}} method.}
    \label{fig:examples-supersonic-n}
\end{figure}

\begin{figure}[!t]
    \centering
    \hspace*{\fill}
    \includegraphics[height=\figureheight]{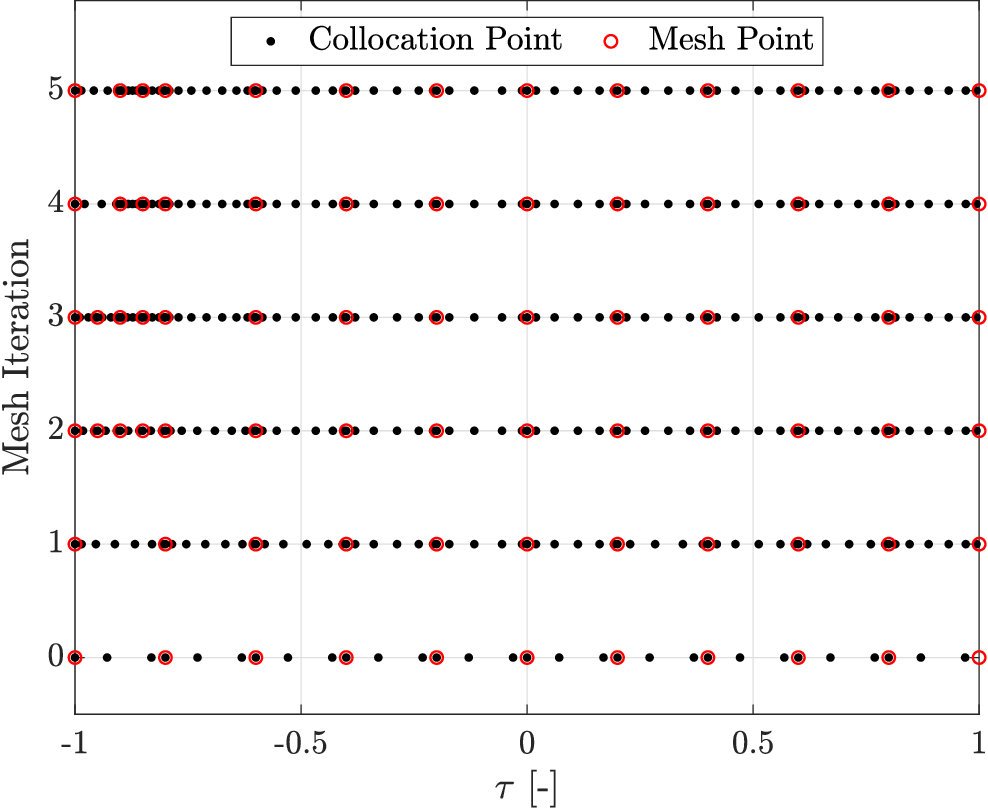}%
    \hspace*{\fill}
    \caption{Mesh refinement history for Example 1 using the \HRvar{3}{10}{\ode{45}} method.}
    \label{fig:examples-supersonic-mesh_1em6}
\end{figure}

Consider the following minimum-time supersonic aircraft climb optimal control problem taken from Ref.~\cite{DarbyRao2011a} using the dynamic model provided in Ref.~\cite{SeywaldWell1994}.
Determine the state, $\mbf{x}(t)=[h(t),v(t),\gamma(t)]$, the control, $u(t)={n}(t)$, and the terminal time, $t_f$, that minimize the objective functional
\begin{equation}\label{eq:examples-supersonic-obj}
    \mcl{J}=t_f,
\end{equation}
subject to the dynamic constraints
\begin{equation}\label{eq:examples-supersonic-dyn}
    \md{h}=v\sin(\gamma),
    \qquad
    \md{v}=\dfrac{T-D}{m}-g\sin(\gamma),
    \qquad
    \md{\gamma}=g\left(\dfrac{n-\cos(\gamma)}{v}\right)\!,
\end{equation}
the boundary conditions
\begin{equation}\label{eq:examples-supersonic-bnd}
    \begin{array}{rcl@{\qquad}rcl@{\qquad}rcl}
        h(0)   & = & 0,      & v(0)   & = & 0.12931, & \gamma(0)   & = & 0, \\[1.0em]
        h(t_f) & = & 19.995, & v(t_f) & = & 0.29509, & \gamma(t_f) & = & 0, \\
    \end{array}
\end{equation}
and the inequality path constraint
\begin{equation}\label{eq:examples-supersonic-pth}
    0\leq\gamma\leq\dfrac{\pi}{2},
\end{equation}
where $h$ is the altitude, $v$ is the speed, $\gamma$ is the flight path angle, $u$ is the load factor, $T$ is the thrust force, $D$ is the drag force, $m$ is the vehicle mass, and $g$ is the acceleration due to gravity.
A solution to this optimal control problem on $\tau\in[-1,+1]$ using the \HRvar{3}{10}{\ode{45}} method is shown in Figs.~\ref{fig:examples-supersonic-state} and~\ref{fig:examples-supersonic-n}.
It is noted that Figs.~\ref{fig:examples-supersonic-h}--\ref{fig:examples-supersonic-gamma} show the optimal state solutions obtained by solving the NLP of Eqs.~\eqref{eq:lgr-obj}--\eqref{eq:lgr-cnt} as well as the solutions obtained via explicit simulation using the \MATLAB~ODE solver \ode{45}, whereas Fig.~\ref{fig:examples-supersonic-n} shows the optimal control solution along with the interpolated control used in the explicit simulation.
% --------------------------------------------------------------------------------------------
% Figure: supersonic example optimal state components vs tau using the phs-(3,10,ode45) method
% --------------------------------------------------------------------------------------------
% ---------------------------------------------------------------------------------------------
% Figure: supersonic example optimal control component vs tau using the phs-(3,10,ode45) method
% ---------------------------------------------------------------------------------------------
For this solution, an optimal final time of $t_f^*=170.565775$ is obtained, and the corresponding mesh refinement history is shown in Fig.~\ref{fig:examples-supersonic-mesh_1em6}, where the \HRvar{3}{10}{\ode{45}} method takes five (5) iterations to satisfy the desired mesh tolerance.
% ------------------------------------------------------------------------------------
% Figure: supersonic example mesh refinement history using the phs-(3,10,ode45) method
% ------------------------------------------------------------------------------------
It is seen in Fig.~\ref{fig:examples-supersonic-gamma} that the state constraint on the flight path angle as given in Eq.~\eqref{eq:examples-supersonic-pth} is active from $\tau=-1$ to $\tau\approx{-0.9}$, and this active constraint arc is also seen in the altitude and load factor solutions of Figs.~\ref{fig:examples-supersonic-h} and~\ref{fig:examples-supersonic-n}, respectively.
To accurately capture this behavior, it is expected that the \HR~method increases the density of the mesh near the activation and deactivation times of the path constraint.
On the initial mesh ($M=0$) of Fig.~\ref{fig:examples-supersonic-mesh_1em6}, the maximum error in the solution is $e_{\max}^{[0]}\approx{3}\times{10}^{-3}>\epsilon$, which is largely due to a purposefully-naive initial mesh.
Within the region of the active state constraint (that is, in $\tau\in[-1,-0.8]$), the first two mesh refinement iterations are driven by the error in the solution of the flight path angle. In this region in Fig.~\ref{fig:examples-supersonic-mesh_1em6}, the first mesh refinement iteration (that is, from mesh $M=0$ to $M=1$) adds collocation points to decrease the error in the solution, and during the second mesh refinement iteration (that is, from mesh $M=1$ to $M=2$), $p$ refinement is exhausted.
As a result, the mesh interval is split.
Then, collocation points are appropriately added until the desired mesh tolerance is satisfied.
For $\tau\in[-0.8,+1]$ on which the solution takes a smooth form (that is, when the state inequality constraint is inactive), the error in the solution is decreased by performing $p$ refinement, and the mesh in $\tau\in[-0.6,+1]$ remains unchanged after the second mesh refinement iteration.

Next, the performance of the \HR~method is compared against the performance of the \PR, \PRHR, \LR, \LRHR, \LRL, and \LRLHR~methods, and the error estimates obtained using the methods of Section~\ref{subsec:mesh-err} and Ref.~\cite{PattersonRao2015} are compared.
% ---------------------------------------------------------------------------
% Table: results for supersonic example using various mesh refinement methods
\begin{table}[!t]
    \centering
    \caption{Results for Example 1 using various mesh refinement methods.}
    \begin{tabular*}{\textwidth}{@{\extracolsep\fill}clccccccc@{}}
        \toprule\toprule
        \multirow{2}{*}{$N_{\min}$} & \multicolumn{1}{c}{\multirow{2}{*}{Method}} & \multirow{2}{*}{$N$} & \multirow{2}{*}{$K$} & \multirow{2}{*}{$M$} & \multirow{2}{*}{$e_{\max}^{[M]}$} & \multirow{2}{*}{$\tilde{e}_{\max}^{[M]}$} & \multicolumn{2}{c}{Total Time [s]} \\
        &&&&&&& SNOPT & CPU \\
        \midrule\midrule
        \multirow{7}{*}{2} & None\textsuperscript{a}          & 20  & 10 & 0 & 1.016$\times{10}^{-1}$ & 9.693$\times{10}^{-2}$ & 0.024 & 0.055 \\
        \cmidrule{2-9}     & \LRvar\textsuperscript{a}        & 107 & 21 & 4 & 1.228$\times{10}^{-6}$ & 9.261$\times{10}^{-7}$ & 0.365 & 0.532 \\
                           & \LRHRvar{\ode{45}}               & 114 & 22 & 4 & 6.195$\times{10}^{-7}$ & 4.696$\times{10}^{-7}$ & 0.411 & 0.978 \\
                           & \LRHRvar{\ode{89}}               & 114 & 22 & 4 & 6.193$\times{10}^{-7}$ & 4.696$\times{10}^{-7}$ & 0.403 & 1.840 \\
        \cmidrule{2-9}     & \HRvar{2}{10}{\ode{45}}          & 83  & 13 & 6 & 9.425$\times{10}^{-7}$ & 7.140$\times{10}^{-7}$ & 0.521 & 1.430 \\
                           & \HRvar{2}{10}{\ode{89}}          & 83  & 13 & 6 & 9.114$\times{10}^{-7}$ & 7.140$\times{10}^{-7}$ & 0.527 & 2.883 \\
        \midrule
        \multirow{7}{*}{3} & None\textsuperscript{a}          & 30  & 10 & 0 & 2.954$\times{10}^{-3}$ & 2.570$\times{10}^{-3}$ & 0.028 & 0.059 \\
        \cmidrule{2-9}     & \LRvar\textsuperscript{a}        & 103 & 16 & 4 & 1.200$\times{10}^{-6}$ & 8.685$\times{10}^{-7}$ & 0.388 & 0.557 \\
                           & \LRHRvar{\ode{45}}               & 99  & 14 & 5 & 5.655$\times{10}^{-7}$ & 4.105$\times{10}^{-7}$ & 0.478 & 1.037 \\
                           & \LRHRvar{\ode{89}}               & 99  & 14 & 5 & 5.655$\times{10}^{-7}$ & 4.105$\times{10}^{-7}$ & 0.474 & 1.809 \\
        \cmidrule{2-9}     & \HRvar{3}{10}{\ode{45}}          & 77  & 12 & 5 & 9.424$\times{10}^{-7}$ & 7.429$\times{10}^{-7}$ & 0.333 & 1.080 \\
                           & \HRvar{3}{10}{\ode{89}}          & 77  & 12 & 5 & 9.462$\times{10}^{-7}$ & 7.429$\times{10}^{-7}$ & 0.349 & 2.381 \\
        \midrule
        \multirow{7}{*}{4} & None\textsuperscript{a}          & 40  & 10 & 0 & 1.269$\times{10}^{-2}$ & 1.010$\times{10}^{-2}$ & 0.039 & 0.070 \\
        \cmidrule{2-9}     & \PRvar{4}{12}\textsuperscript{a} & 82  & 13 & 4 & 9.428$\times{10}^{-7}$ & 7.395$\times{10}^{-7}$ & 0.343 & 0.502 \\
                           & \PRHRvar{4}{12}{\ode{45}}        & 85  & 13 & 4 & 9.429$\times{10}^{-7}$ & 7.395$\times{10}^{-7}$ & 0.357 & 0.790 \\
                           & \PRHRvar{4}{12}{\ode{89}}        & 85  & 13 & 4 & 9.419$\times{10}^{-7}$ & 7.395$\times{10}^{-7}$ & 0.359 & 1.468 \\
        \cmidrule{2-9}     & \HRvar{4}{12}{\ode{45}}          & 77  & 12 & 5 & 9.425$\times{10}^{-7}$ & 7.400$\times{10}^{-7}$ & 0.375 & 1.049 \\
                           & \HRvar{4}{12}{\ode{89}}          & 77  & 12 & 5 & 9.425$\times{10}^{-7}$ & 7.400$\times{10}^{-7}$ & 0.375 & 2.100 \\
        \bottomrule\bottomrule
    \end{tabular*}
    \begin{tablenotes}
        \item[$^{\rm{a}}$] The error estimate $e_{\max}^{[M]}$ of Section~\ref{subsec:mesh-err} is obtained with the \MATLAB~ODE solver \ode{45}.
    \end{tablenotes}
    \label{tab:examples-supersonic-results}
\end{table}%
% ---------------------------------------------------------------------------
Table~\ref{tab:examples-supersonic-results} summarizes the results for Example 1 using various mesh refinement methods, where only the methods that obtain the smallest final mesh (that is, with the fewest total number of collocation points) are shown.
In terms of final mesh size, the \HR~method converges to the smallest mesh for every considered value of $N_{\min}$ when compared against that obtained using previously developed methods, where a 22\%, 22\%, and 6\% decrease in total number of collocation points and a 38\%, 14\%, and 8\% decrease in total number of mesh intervals are observed for $N_{\min}=[2,3,4]$, respectively.
The decrease in mesh size, however, typically requires the \HR~method to perform more mesh refinement iterations to satisfy the desired mesh tolerance.
The \HR~method also requires the use of an explicit simulation scheme, which is typically more computationally expensive than the Gaussian quadrature integration scheme utilized in the \PR, \LR, and \LRL~methods.
As a result, expected increases in total CPU times are observed when using the \HR~method.
In terms of error analysis, the error estimate $e_{\max}^{[M]}$ obtained in Section~\ref{subsec:mesh-err} is larger than the error estimate $\tilde{e}_{\max}^{[M]}$ obtained via the method of Ref.~\cite{PattersonRao2015} for all results shown in Table~\ref{tab:examples-supersonic-results}.
For $N_{\min}=[2,3]$, the final meshes obtained using the \LR~method do not satisfy $e_{\max}^{[M]}\leq\epsilon$, which suggests there is still a discrepancy between solutions obtained via collocation and time-marching schemes on these final meshes.
Similar phenomena are observed in many other cases for Example 1 when using the \PR, \LR, and \LRL~methods but are not shown in Table~\ref{tab:examples-supersonic-results}.
On the other hand, all final meshes obtained using the \HR~method satisfy both $e_{\max}^{[M]}\leq\epsilon$ and $\tilde{e}_{\max}^{[M]}\leq\epsilon$, where the solutions obtained on the final mesh via collocation and time-marching are guaranteed to be in agreement with one another.
Finally, the choice of \MATLAB~ODE solver did not impact the size of the final mesh for this example for all values of $N_{\min}$ considered, as shown in Table~\ref{tab:examples-supersonic-results}; however, using the \MATLAB~ODE solver \ode{89} required larger computational times than \ode{45}.
To summarize for Example 1, the \HR~mesh refinement method is able to outperform the \PR, \PRHR, \LR, \LRHR, \LRL, and \LRLHR~mesh refinement methods in terms of final mesh size.

% ------------------------------------------------------------------------------------------------ %
% ---------- SUBSECTION: Example 1: Minimum-Time to Climb of a Supersonic Aircraft (END) --------- %
% ------------------------------------------------------------------------------------------------ %

% ------------------------------------------------------------------------------------------------ %
% ----------- SUBSECTION: Example 2: Minimum-Time Reorientation of a Robot Arm (START) ----------- %
% ------------------------------------------------------------------------------------------------ %
\subsection{Example 2: Minimum-Time Reorientation of a Robot Arm}\label{subsec:examples-robotarm}

\begin{figure}[!t]
    \centering
    \hspace*{\fill}
    \subfloat[$\rho(t)$ vs. $\tau$.]{%
        \includegraphics[height=\figureheight]{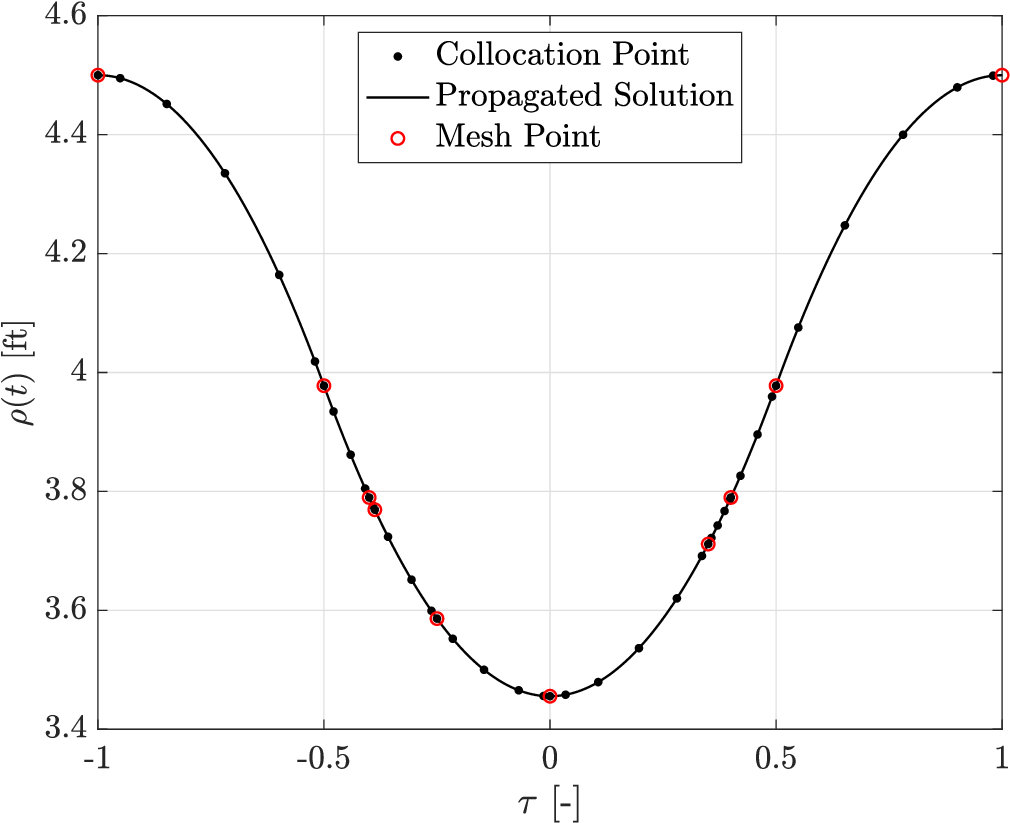}%
        \label{fig:examples-robotarm-rho}
        }
    \hspace*{\fill}
    \subfloat[$\theta(t)$ vs. $\tau$.]{%
        \includegraphics[height=\figureheight]{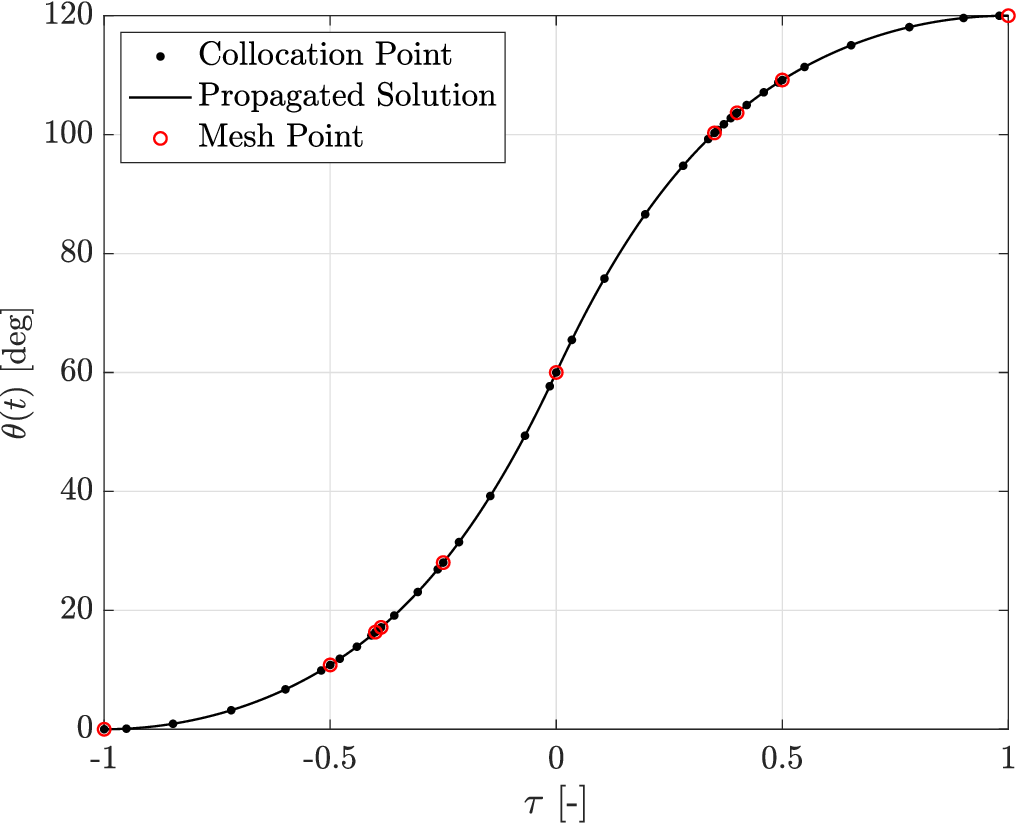}%
        \label{fig:examples-robotarm-theta}
        }
    \hspace*{\fill}
    \\
    \hspace*{\fill}
    \subfloat[$\phi(t)$ vs. $\tau$.]{%
        \includegraphics[height=\figureheight]{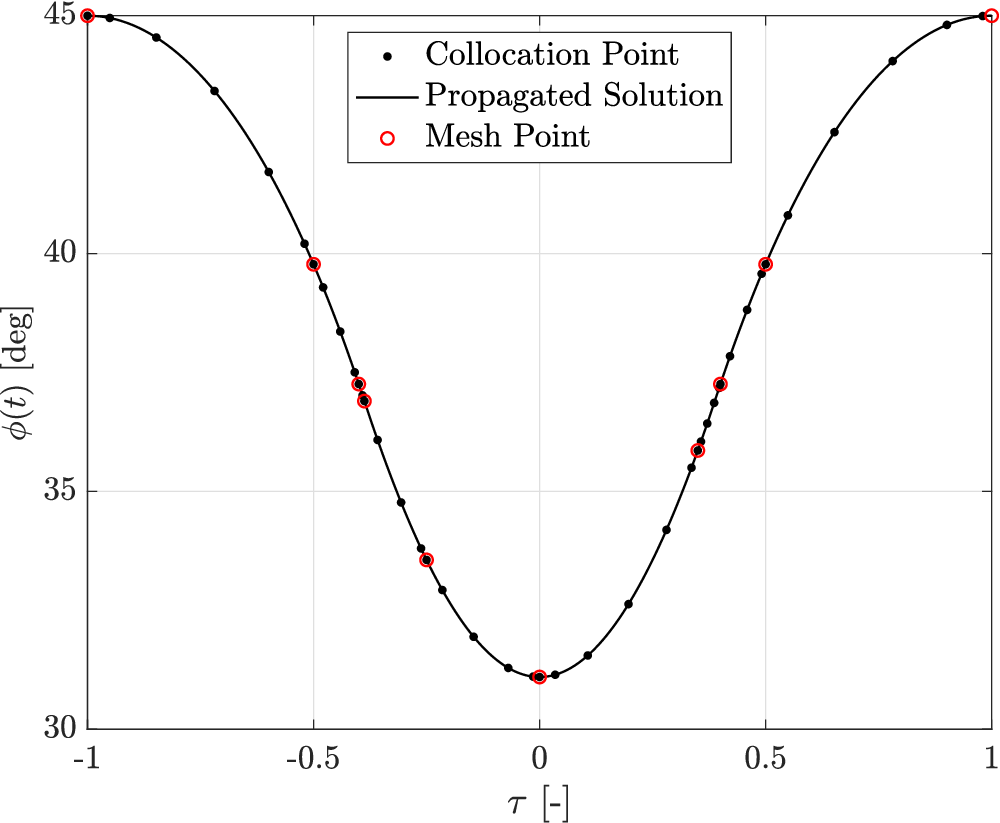}%
        \label{fig:examples-robotarm-phi}
        }
    \hspace*{\fill}
    \subfloat[$\md{\rho}(t)$ vs. $\tau$.]{%
        \includegraphics[height=\figureheight]{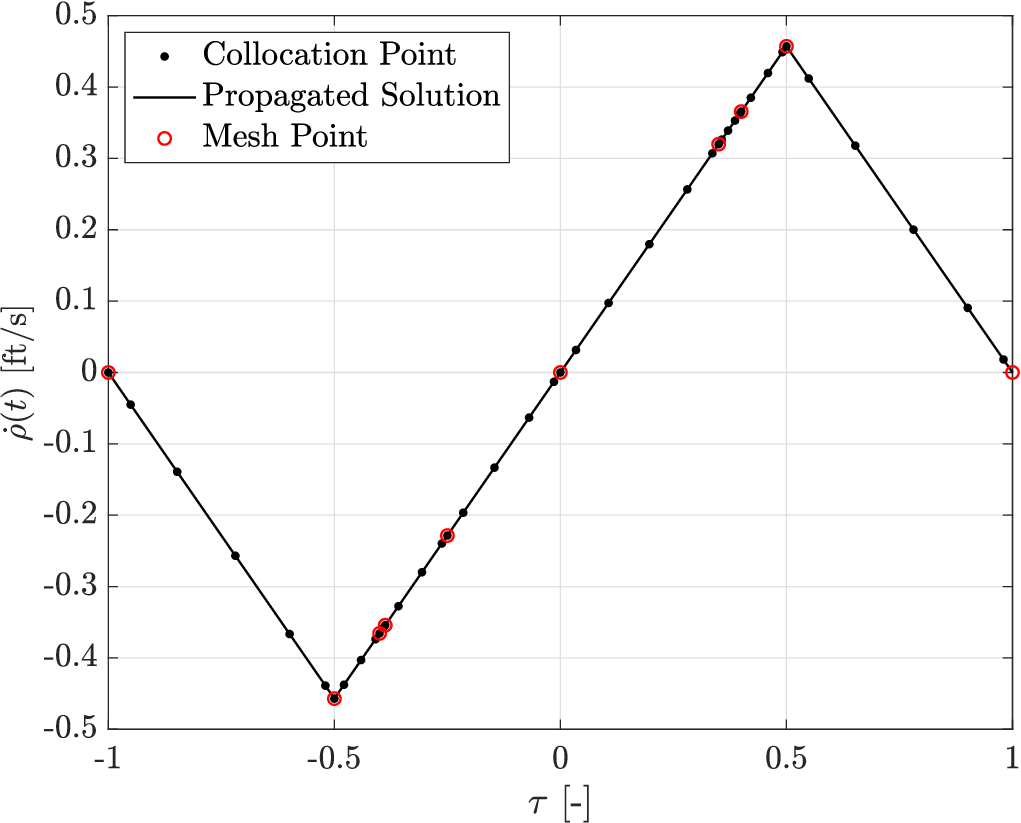}%
        \label{fig:examples-robotarm-rhodot}
        }
    \hspace*{\fill}
    \\
    \hspace*{\fill}
    \subfloat[$\md{\theta}(t)$ vs. $\tau$.]{%
        \includegraphics[height=\figureheight]{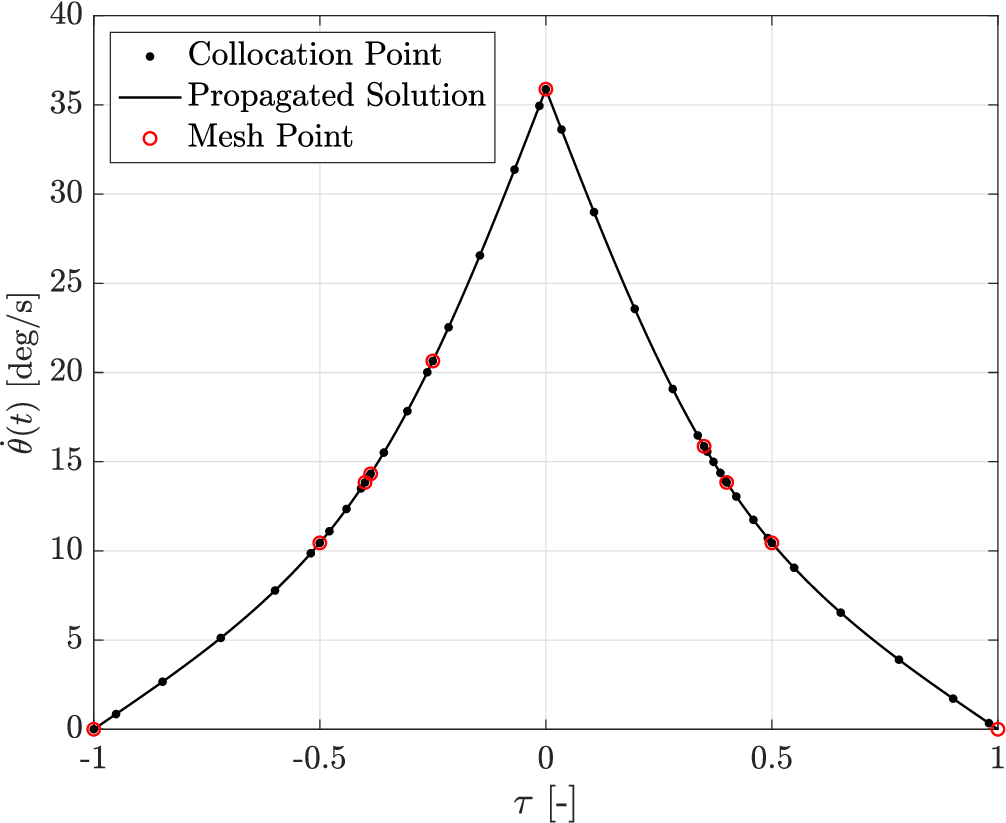}%
        \label{fig:examples-robotarm-thetadot}
        }
    \hspace*{\fill}
    \subfloat[$\md{\phi}(t)$ vs. $\tau$.]{%
        \includegraphics[height=\figureheight]{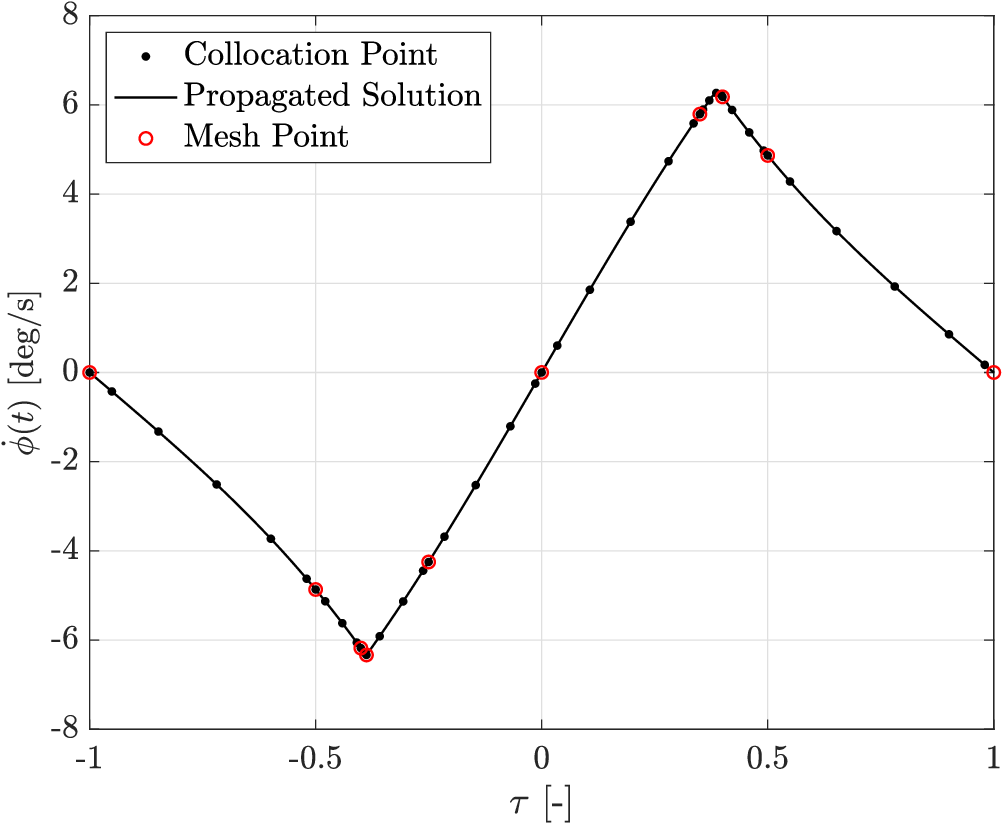}%
        \label{fig:examples-robotarm-phidot}
        }
    \hspace*{\fill}
    \caption{Optimal state components for Example 2 using the \HRvar{2}{6}{\ode{45}} method.}
    \label{fig:examples-robotarm-state}
\end{figure}

\begin{figure}[!t]
    \centering
    \hspace*{\fill}
    \subfloat[$u_{\rho}(t)$ vs. $\tau$.]{%
        \includegraphics[height=\figureheight]{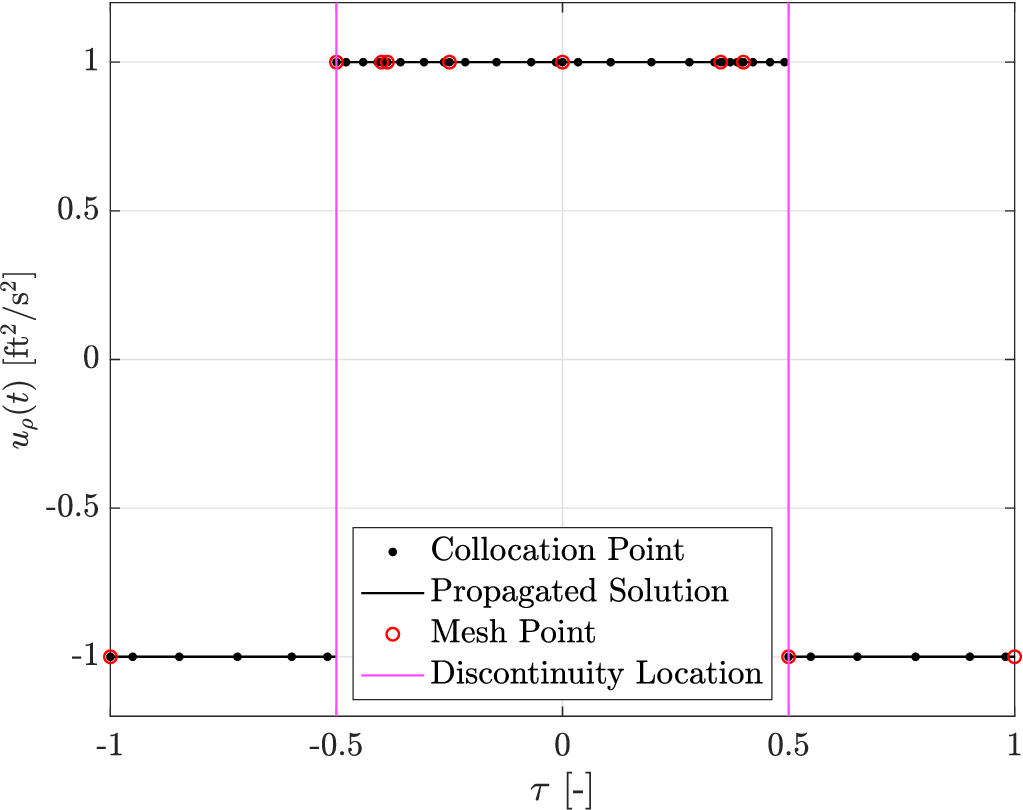}%
        \label{fig:examples-robotarm-urho}
        }
    \hspace*{\fill}
    \subfloat[$u_{\theta}(t)$ vs. $\tau$.]{%
        \includegraphics[height=\figureheight]{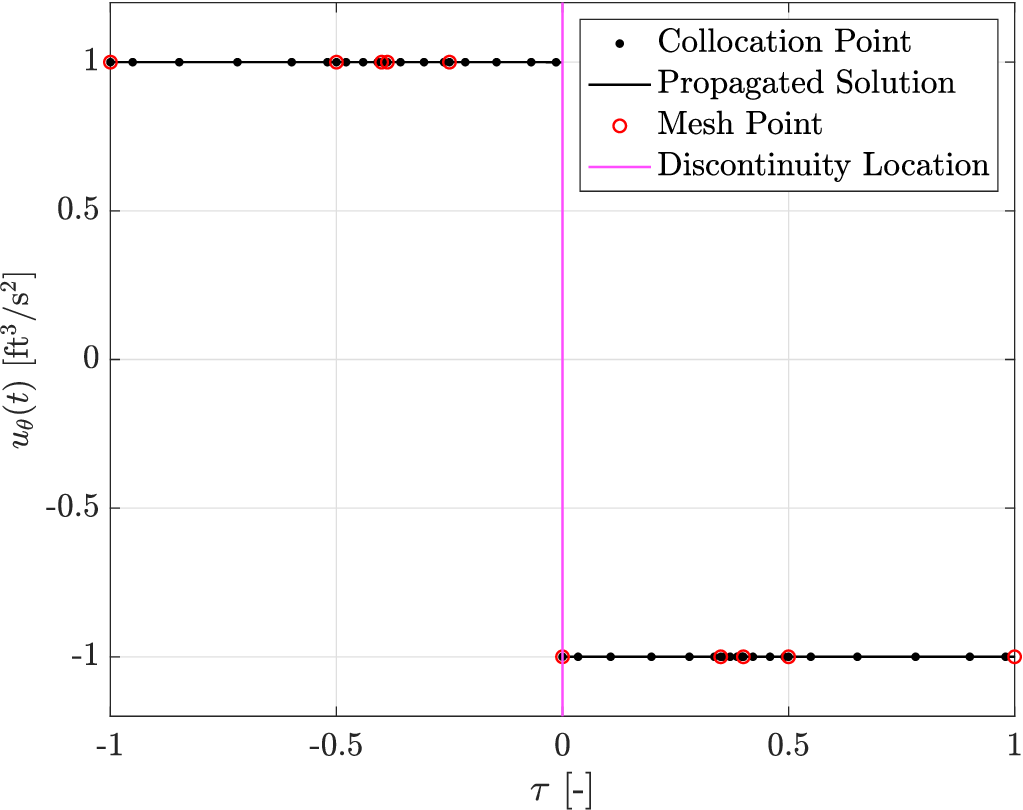}%
        \label{fig:examples-robotarm-utheta}
        }
    \hspace*{\fill}
    \\
    \hspace*{\fill}
    \subfloat[$u_{\phi}(t)$ vs. $\tau$.]{%
        \includegraphics[height=\figureheight]{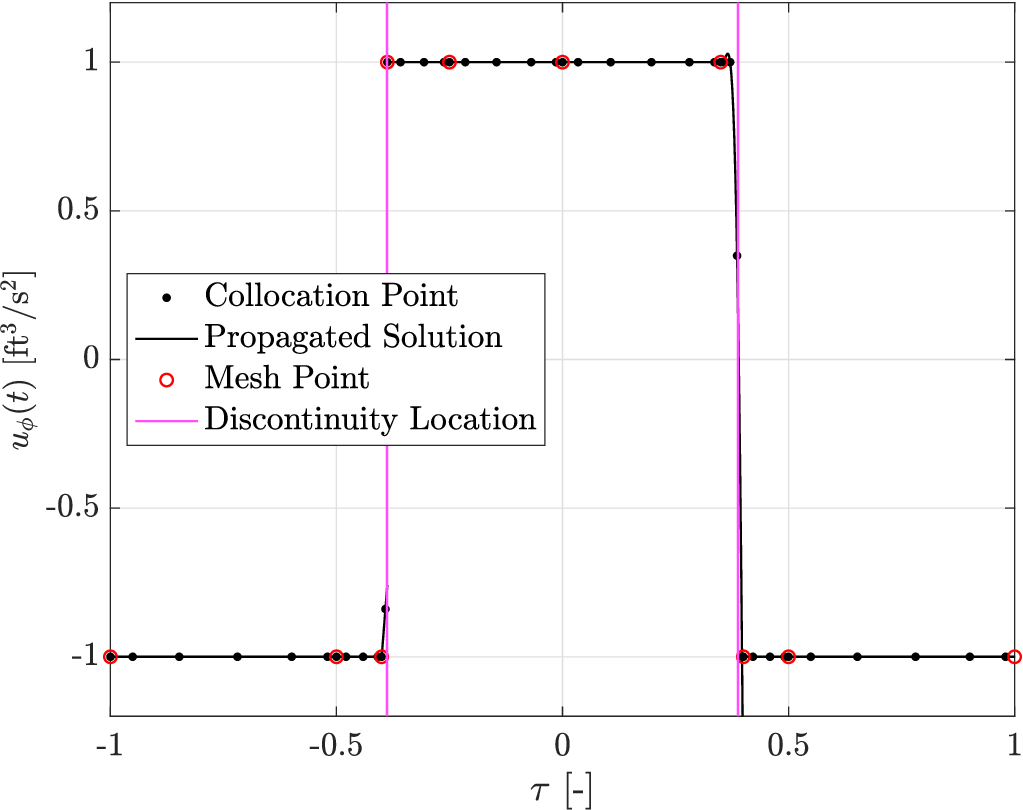}%
        \label{fig:examples-robotarm-uphi}
        }
    \hspace*{\fill}
    \caption{Optimal control components for Example 2 using the \HRvar{2}{6}{\ode{45}} method.}
    \label{fig:examples-robotarm-control}
\end{figure}

\begin{figure}[!t]
    \centering
    \hspace*{\fill}
    \includegraphics[height=\figureheight]{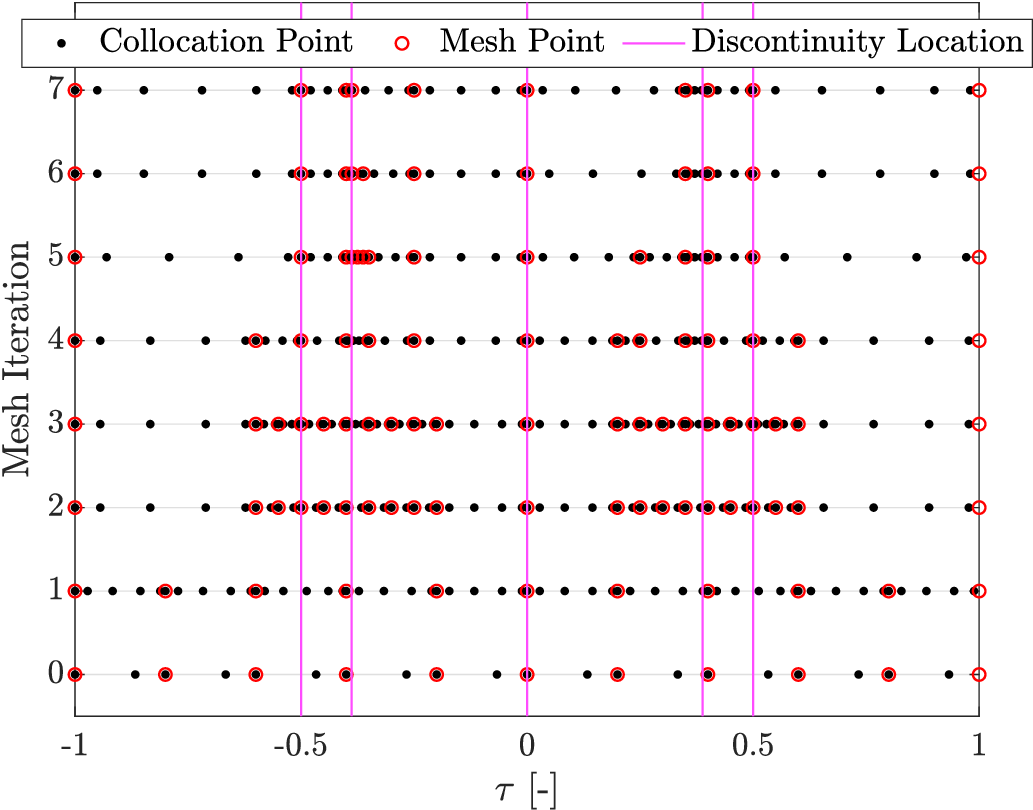}%
    \hspace*{\fill}
    \caption{Mesh refinement history for Example 2 using the \HRvar{2}{6}{\ode{45}} method.}
    \label{fig:examples-robotarm-mesh_1em6}
\end{figure}

Consider the following minimum-time robot arm reorientation optimal control problem taken from Ref.~\cite{DolanMunson2004}.
Determine the state, $\mbf{x}(t)=[\rho(t),\theta(t),\phi(t),\md{\rho}(t),\md{\theta}(t),\md{\phi}(t)]$, the control, $\mbf{u}(t)=[u_{\rho}(t),u_{\theta}(t),u_{\phi}(t)]$, and the terminal time, $t_f$, that minimize the objective functional
\begin{equation}\label{eq:examples-robotarm-obj}
    \mcl{J}=t_f,
\end{equation}
subject to the dynamic constraints
\begin{equation}\label{eq:examples-robotarm-dyn}
    \mdd{\rho}=\dfrac{u_{\rho}}{L},
    \qquad
    \mdd{\theta}=\dfrac{u_{\theta}}{I_{\phi}\sin^2(\phi)},
    \qquad
    \mdd{\phi}=\dfrac{u_{\phi}}{I_{\phi}},
\end{equation}
the boundary conditions
\begin{equation}\label{eq:examples-robotarm-bnd}
    \begin{array}{rcl@{\qquad}rcl@{\qquad}rcl@{\qquad}rcl@{\qquad}rcl@{\qquad}rcl}
        \rho(0)   & = & \dfrac{9}{2}, & \theta(0)   & = & 0,               & \phi(0)   & = & \dfrac{\pi}{4}, & \md{\rho}(0)   & = & 0, & \md{\theta}(0)   & = & 0, & \md{\phi}(0)   & = & 0, \\[1.0em]
        \rho(t_f) & = & \dfrac{9}{2}, & \theta(t_f) & = & \dfrac{2\pi}{3}, & \phi(t_f) & = & \dfrac{\pi}{4}, & \md{\rho}(t_f) & = & 0, & \md{\theta}(t_f) & = & 0, & \md{\phi}(t_f) & = & 0, \\
    \end{array}
\end{equation}
and the inequality path constraints
\begin{equation}\label{eq:examples-robotarm-pth}
    \begin{array}{rcccl@{\qquad}rcccl@{\qquad}rcccl}
        0  & \leq & \rho     & \leq & L, & -\pi & \leq & \theta     & \leq & \pi, & 0  & \leq & \phi     & \leq & \pi, \\[1.0em]
        -1 & \leq & u_{\rho} & \leq & 1, & -1   & \leq & u_{\theta} & \leq & 1,   & -1 & \leq & u_{\phi} & \leq & 1, \\
    \end{array}
\end{equation}
where $L=5$ and $I_{\phi}\equiv{I}_{\phi}(\rho)=((L-\rho)^3+\rho^3)/3$.
A solution to this optimal control problem on $\tau\in[-1,+1]$ using the \HRvar{2}{6}{\ode{45}} method is shown in Figs.~\ref{fig:examples-robotarm-state} and~\ref{fig:examples-robotarm-control}.
It is noted that Figs.~\ref{fig:examples-robotarm-rho}--\ref{fig:examples-robotarm-phidot} show the optimal state solutions obtained by solving the NLP of Eqs.~\eqref{eq:lgr-obj}--\eqref{eq:lgr-cnt} as well as the solutions obtained via explicit simulation using the \MATLAB~ODE solver \ode{45}, whereas Figs.~\ref{fig:examples-robotarm-urho}--\ref{fig:examples-robotarm-uphi} show the optimal control solutions along with the interpolated controls used in the explicit simulation.
% -----------------------------------------------------------------------------------------
% Figure: robot arm example optimal state component vs tau using the phs-(2,6,ode45) method
% -----------------------------------------------------------------------------------------
% -------------------------------------------------------------------------------------------
% Figure: robot arm example optimal control component vs tau using the phs-(2,6,ode45) method
% -------------------------------------------------------------------------------------------
For this solution, an optimal final time of $t_f^{*}=9.140963$ is obtained, and the corresponding mesh refinement history is shown in Fig.~\ref{fig:examples-robotarm-mesh_1em6}, where the \HRvar{2}{6}{\ode{45}} method takes seven (7) iterations to satisfy the desired mesh tolerance.
% ----------------------------------------------------------------------------------
% Figure: robot arm example mesh refinement history using the phs-(2,6,ode45) method
% ----------------------------------------------------------------------------------
It is seen in Fig.~\ref{fig:examples-robotarm-control} that the optimal control components all take on a bang-bang structure, where it is known that the five control discontinuities are located at $\tau\approx\{-0.5,-0.3882,0,+0.3882,+0.5\}$.
These discontinuity locations also correspond to the local extrema in the respective state components $\md{\rho}$, $\md{\theta}$, and $\md{\phi}$ shown in Figs.~\ref{fig:examples-robotarm-rhodot}--\ref{fig:examples-robotarm-phidot}, respectively.
Because of the rapid changes in the solution, it is expected that the \HR~method increases the density of the mesh near the locations of the control discontinuities.
On the initial mesh ($M=0$) of Fig.~\ref{fig:examples-robotarm-mesh_1em6}, the maximum error in the solution is $e_{\max}^{[0]}\approx{1.5}\times{10}^{-3}>\epsilon$, which is largely due to an initial mesh that does not accurately capture all of the discontinuity locations; however, due to the construction of the initial mesh consisting of 10 uniformly spaced intervals, the discontinuity in the control component $u_{\theta}$ at $\tau=0$ is already identified.
Mesh refinement is initially driven by the error in the solutions of $\phi$ and $\md{\phi}$; however, $p$ refinement is quickly exhausted in $\tau\in[-0.6,-0.2]$ and $\tau\in[+0.2,+0.6]$ during the second mesh refinement iteration (that is, from mesh $M=1$ to $M=2$) due to large state error estimates from the remaining inaccurately captured discontinuities.
As a result, the appropriate intervals are split, and the discontinuities in the control component $u_{\rho}$ at $\tau=\{-0.5,+0.5\}$ are then identified.
It is important to note that the $h$ refinement step of Section~\ref{subsubsec:mesh-refine-hp} does not perform any discontinuity detection; thus, new intervals may be formed in which discontinuities are not present.
The $h$ reduction step of Section~\ref{subsubsec:mesh-refine-hm} then helps remove unnecessary mesh intervals during the fourth and fifth mesh refinement iterations (that is, from mesh $M=3$ to $M=4$ and from $M=4$ to $M=5$, respectively) and, as a result, keeps the size of the mesh small.
Finally, appropriate regions of the solution are further refined until the remaining two discontinuities in the control component $u_{\phi}$ at $\tau\approx\{-0.3882,+0.3882\}$ are accurately identified, and the desired mesh tolerance is satisfied.

Next, the performance of the various methods and error estimates is compared in a manner similar to that done for Example 1 in Section~\ref{subsec:examples-supersonic}.
Table~\ref{tab:examples-robotarm-results} summarizes the results for Example 2, where only the methods that obtain the smallest final mesh are shown.
% --------------------------------------------------------------------------
% Table: results for robot arm example using various mesh refinement methods
\begin{table}[!t]
    \centering
    \caption{Results for Example 2 using various mesh refinement methods.}
    \begin{tabular*}{\textwidth}{@{\extracolsep\fill}clccccccc@{}}
        \toprule\toprule
        \multirow{2}{*}{$N_{\min}$} & \multicolumn{1}{c}{\multirow{2}{*}{Method}} & \multirow{2}{*}{$N$} & \multirow{2}{*}{$K$} & \multirow{2}{*}{$M$} & \multirow{2}{*}{$e_{\max}^{[M]}$} & \multirow{2}{*}{$\tilde{e}_{\max}^{[M]}$} & \multicolumn{2}{c}{Total Time [s]} \\
        &&&&&&& SNOPT & CPU \\
        \midrule\midrule
        \multirow{7}{*}{2} & None\textsuperscript{a}          & 20 & 10 & 0 & 1.505$\times{10}^{-3}$ & 1.125$\times{10}^{-3}$ & 0.018 & 0.049 \\
        \cmidrule{2-9}     & \LRvar\textsuperscript{a}        & 84 & 17 & 4 & 2.479$\times{10}^{-7}$ & 1.894$\times{10}^{-7}$ & 0.140 & 0.307 \\
                           & \LRHRvar{\ode{45}}               & 84 & 17 & 4 & 2.479$\times{10}^{-7}$ & 1.894$\times{10}^{-7}$ & 0.139 & 0.508 \\
                           & \LRHRvar{\ode{89}}               & 84 & 17 & 4 & 2.468$\times{10}^{-7}$ & 1.894$\times{10}^{-7}$ & 0.138 & 0.992 \\
        \cmidrule{2-9}     & \HRvar{2}{6}{\ode{45}}           & 42 & 9  & 7 & 6.615$\times{10}^{-7}$ & 5.197$\times{10}^{-7}$ & 2.806 & 3.432 \\
                           & \HRvar{2}{6}{\ode{89}}           & 42 & 9  & 7 & 6.671$\times{10}^{-7}$ & 5.197$\times{10}^{-7}$ & 2.823 & 4.310 \\
        \midrule
        \multirow{7}{*}{3} & None\textsuperscript{a}          & 30 & 10 & 0 & 1.902$\times{10}^{-4}$ & 1.690$\times{10}^{-4}$ & 0.020 & 0.052 \\
        \cmidrule{2-9}     & \PRvar{3}{10}\textsuperscript{a} & 81 & 19 & 6 & 2.250$\times{10}^{-6}$ & 3.498$\times{10}^{-7}$ & 0.174 & 0.398 \\
                           & \PRHRvar{3}{8}{\ode{45}}         & 84 & 20 & 4 & 4.810$\times{10}^{-7}$ & 3.770$\times{10}^{-7}$ & 0.108 & 0.502 \\
                           & \PRHRvar{3}{8}{\ode{89}}         & 84 & 20 & 4 & 4.808$\times{10}^{-7}$ & 3.770$\times{10}^{-7}$ & 0.109 & 1.081 \\
        \cmidrule{2-9}     & \HRvar{3}{10}{\ode{45}}          & 47 & 9  & 9 & 6.615$\times{10}^{-7}$ & 5.197$\times{10}^{-7}$ & 0.166 & 0.989 \\
                           & \HRvar{3}{10}{\ode{89}}          & 57 & 11 & 7 & 9.017$\times{10}^{-7}$ & 6.845$\times{10}^{-7}$ & 0.144 & 1.449 \\
        \midrule
        \multirow{7}{*}{4} & None\textsuperscript{a}          & 40 & 10 & 0 & 1.661$\times{10}^{-4}$ & 1.343$\times{10}^{-4}$ & 0.024 & 0.057 \\
        \cmidrule{2-9}     & \PRvar{4}{6}\textsuperscript{a}  & 71 & 17 & 6 & 1.193$\times{10}^{-6}$ & 9.174$\times{10}^{-7}$ & 0.111 & 0.332 \\
                           & \PRHRvar{4}{6}{\ode{45}}         & 72 & 17 & 7 & 7.536$\times{10}^{-7}$ & 6.368$\times{10}^{-7}$ & 0.129 & 0.746 \\
                           & \PRHRvar{4}{6}{\ode{89}}         & 72 & 17 & 7 & 7.555$\times{10}^{-7}$ & 6.368$\times{10}^{-7}$ & 0.130 & 1.595 \\
        \cmidrule{2-9}     & \HRvar{4}{6}{\ode{45}}           & 44 & 9  & 5 & 9.420$\times{10}^{-7}$ & 7.175$\times{10}^{-7}$ & 0.090 & 0.525 \\
                           & \HRvar{4}{6}{\ode{89}}           & 44 & 9  & 5 & 9.420$\times{10}^{-7}$ & 7.175$\times{10}^{-7}$ & 0.092 & 1.033 \\
        \bottomrule\bottomrule
    \end{tabular*}
    \begin{tablenotes}
        \item[$^{\rm{a}}$] The error estimate $e_{\max}^{[M]}$ of Section~\ref{subsec:mesh-err} is obtained with the \MATLAB~ODE solver \ode{45}.
    \end{tablenotes}
    \label{tab:examples-robotarm-results}
\end{table}%
% --------------------------------------------------------------------------
In terms of final mesh size, the \HR~method converges to the smallest mesh for every considered value of $N_{\min}$ when compared against previously developed methods, where a 50\%, 30\%, and 38\% decrease in total number of collocation points and a 47\%, 42\%, and 47\% decrease in total number of mesh intervals are observed for $N_{\min}=[2,3,4]$, respectively.
Again, the decrease in mesh size typically requires the \HR~method to perform more mesh refinement iterations to satisfy the desired mesh tolerance; however, the case for $N_{\min}=4$ demonstrates that the \HR~method has the potential to also outperform previously developed methods in terms of computational efficiency.
This example also reiterates trends observed in Section~\ref{subsec:examples-supersonic}, specifically where the \HR~method satisfies $e_{\max}^{[M]}\leq\epsilon$ and $\tilde{e}_{\max}^{[M]}\leq\epsilon$ for every case shown in Table~\ref{tab:examples-robotarm-results}, but the \PR, \LR, and \LRL~methods are unable to guarantee satisfaction of $e_{\max}^{[M]}\leq\epsilon$.
For $N_{\min}=3$, performance of the \HR~method is impacted by the choice of \MATLAB~ODE solver, where the use of \ode{89} requires a slightly larger final mesh compared to \ode{45}.
To summarize for Example 2, the \HR~mesh refinement method again demonstrates a significant improvement over the \PR, \PRHR, \LR, \LRHR, \LRL, and \LRLHR~mesh refinement methods in terms of obtaining a smaller final mesh that satisfies the desired mesh tolerance.

% ------------------------------------------------------------------------------------------------ %
% ------------ SUBSECTION: Example 2: Minimum-Time Reorientation of a Robot Arm (END) ------------ %
% ------------------------------------------------------------------------------------------------ %

% ------------------------------------------------------------------------------------------------ %
% -------------------- SUBSECTION: Example 3: Hyper-Sensitive Problem (START) -------------------- %
% ------------------------------------------------------------------------------------------------ %
\subsection{Example 3: Hyper-Sensitive Problem}\label{subsec:examples-hypersensitive}

Consider the following variation of the hyper-sensitive optimal control problem taken from Refs.~\cite{RaoMease2000,DarbyRao2011a}.
Determine the state, $x(t)$, and the control, $u(t)$, that minimize the objective functional
\begin{equation}\label{eq:examples-hypersensitive-obj}
    \mcl{J}=\dfrac{1}{2}\int_{0}^{t_f}(x^2+u^2)\dt,
\end{equation}
subject to the dynamic constraints
\begin{equation}\label{eq:examples-hypersensitive-dyn}
    \md{x}=-x^3+u,
\end{equation}
and the boundary conditions
\begin{equation}\label{eq:examples-hypersensitive-bnd}
    \begin{array}{rcl@{\qquad}rcl}
        x(0) & = & 1.5, & x(t_f) & = & 1, \\
    \end{array}
\end{equation}
where a final time of $t_f=10,000$ is used.
For sufficiently large values of $t_f$, the resulting Hamiltonian boundary-value problem (HBVP) (that is, the HBVP obtained from deriving the first-order necessary conditions for optimality) is completely {\em{hyper-sensitive}} and the solution exhibits a {\em{take-off}}, {\em{cruise}}, and {\em{landing}} structure~\cite{RaoMease2000}.
It is noted that the control is essentially zero during the cruise segment because the cruise segment lies in the neighborhood of an equilibrium point.

It is known that numerical difficulties often arise when obtaining solutions to the optimal control problem of Eqs.~\eqref{eq:examples-hypersensitive-obj}--\eqref{eq:examples-hypersensitive-bnd} using an explicit simulation method (for example, indirect or direct shooting) due to instabilities in the dynamics in either forward or backward time.
Because the error estimate of Section~\ref{subsec:mesh-err} employs explicit simulation, it, too, may be prone to  numerical difficulties similar to that of indirect or direct shooting.
The objective of this example is to highlight the impact of the initial mesh and ODE solver on the ability to employ the \HR~method due to the reliance of the error estimate on explicit simulation. 

To demonstrate the potential issues with the reliance of the error estimate of the \HR~method on explicit simulation, consider the explicit simulation of the uncontrolled dynamics $\md{x}=-x^3$.
Using any terminal condition perturbed significantly from $x=0$ to solve the TVP of Eq.~\eqref{eq:mesh-sim-tvp}, every \MATLAB~ODE solver fails in backward time because $\md{x}=-x^3$ is unstable in backward time.
As a workaround to the inability to perform explicit simulation in backward time on this example, the contribution to the maximum relative error estimate in Eq.~\eqref{eq:mesh-err-relmax} from the single-interval TVP of Eq.~\eqref{eq:mesh-sim-tvp} is ignored in the mesh refinement method of Algorithm~\ref{alg:mesh-sum} for this example.
A similar approach is applied to the $h$ reduction step of Section~\ref{subsubsec:mesh-refine-hm}, where the contribution to the maximum relative error estimate in Eq.~\eqref{eq:mesh-refine-hm-relmax} from the multiple-interval TVP of Eq.~\eqref{eq:mesh-refine-hm-tvp} is ignored.
While it is not possible to include the contribution to the maximum relative error estimate from the solution of the TVP for this example, accurate error estimates are still obtained in each mesh interval by solving the single- and multiple-interval IVPs of Eqs.~\eqref{eq:mesh-sim-ivp} and~\eqref{eq:mesh-refine-hm-ivp}.
Finally, it is noted that same approaches are utilized for this example when obtaining the error estimates using the \PRHR, \LRHR, and \LRLHR~methods.

\begin{figure}[!t]
    \centering
    \hspace*{\fill}
    \subfloat[Near $\tau=-1$.]{%
        \includegraphics[height=\figureheight]{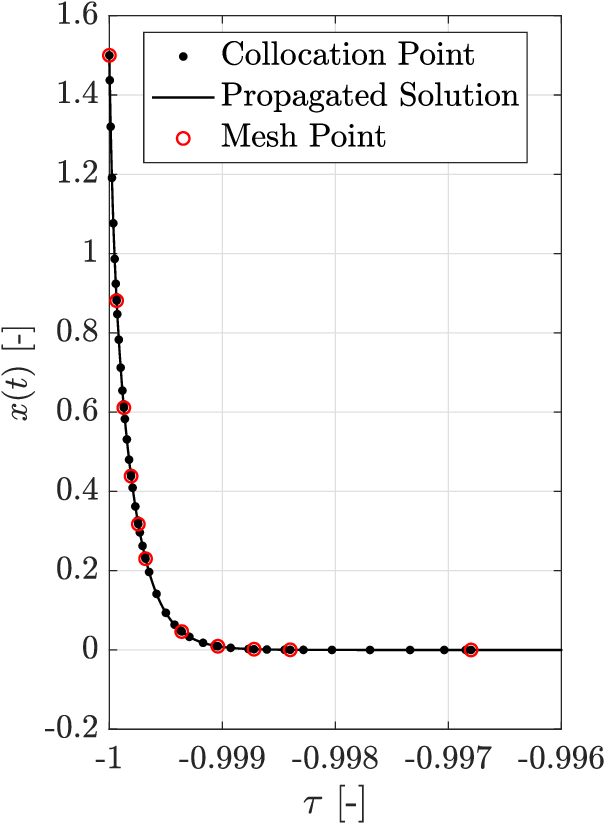}%
        \label{fig:examples-hypersensitive-xstart}
        }
    \hspace*{\fill}
    \subfloat[Over $\tau\in\ensuremath{[-1,+1]}$.]{%
        \includegraphics[height=\figureheight]{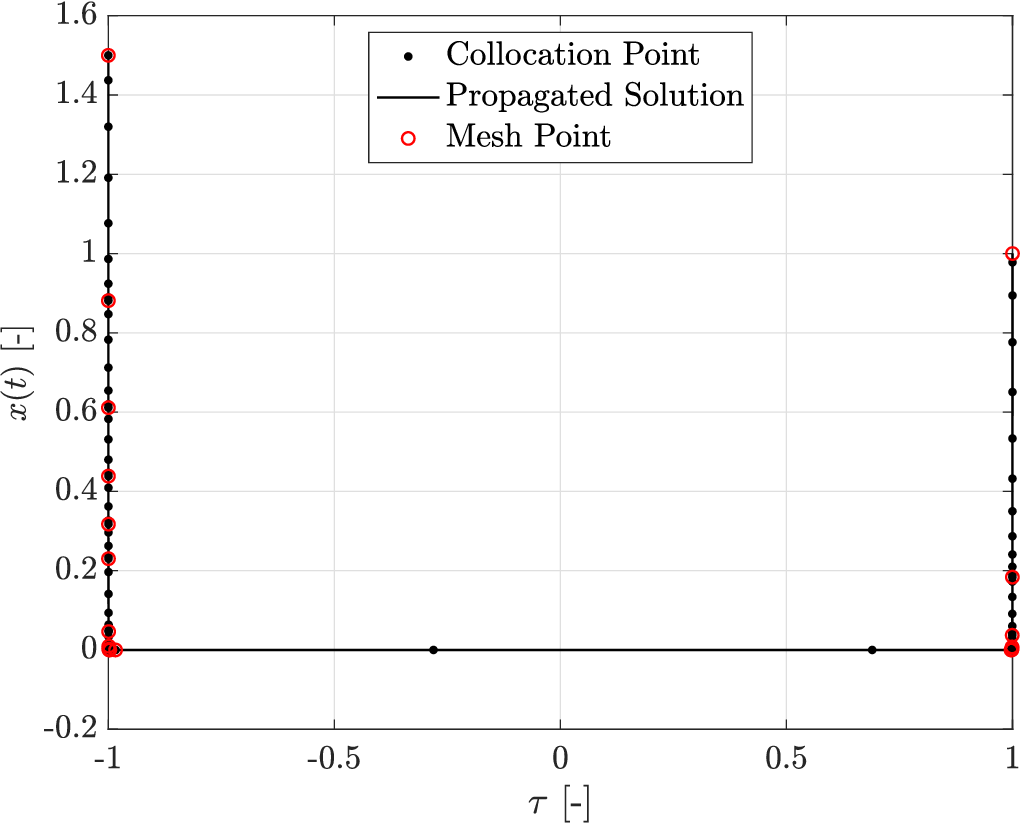}%
        \label{fig:examples-hypersensitive-x}
        }
    \hspace*{\fill}
    \subfloat[Near $\tau=+1$.]{%
        \includegraphics[height=\figureheight]{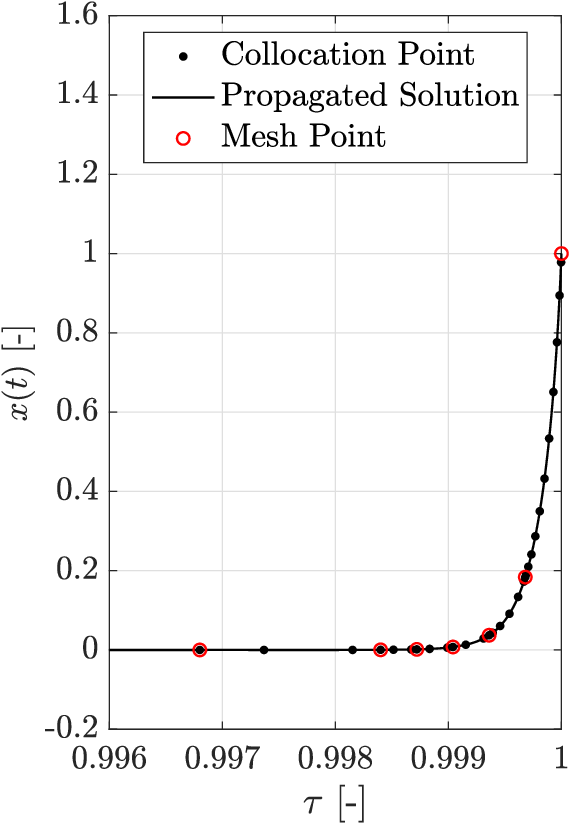}%
        \label{fig:examples-hypersensitive-xend}
        }
    \hspace*{\fill}
    \caption{Optimal state, $x(t)$ vs. $\tau$, for Example 3 using the \HRvar{3}{12}{\ode{89}} method.}
    \label{fig:examples-hypersensitive-state}
\end{figure}

\begin{figure}[!t]
    \centering
    \hspace*{\fill}
    \subfloat[Near $\tau=-1$.]{%
        \includegraphics[height=\figureheight]{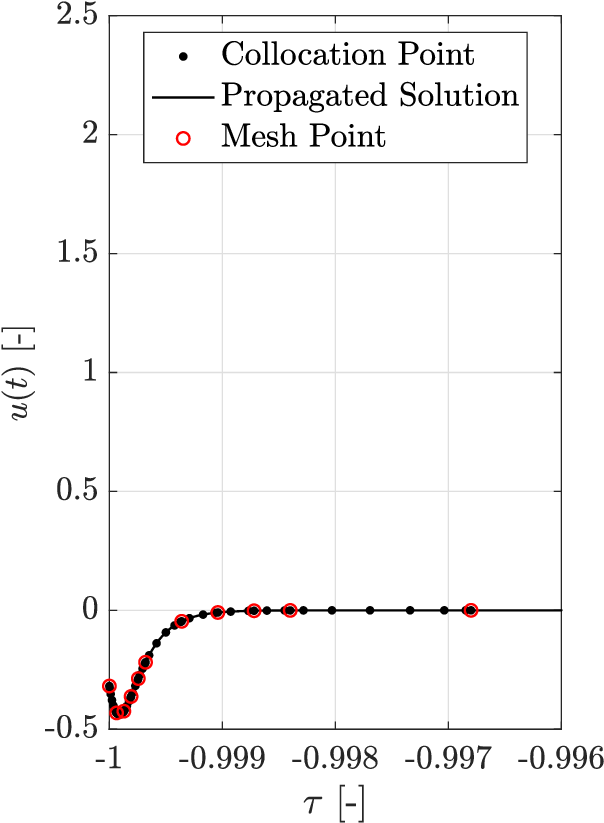}%
        \label{fig:examples-hypersensitive-ustart}
        }
    \hspace*{\fill}
    \subfloat[Over $\tau\in\ensuremath{[-1,+1]}$.]{%
        \includegraphics[height=\figureheight]{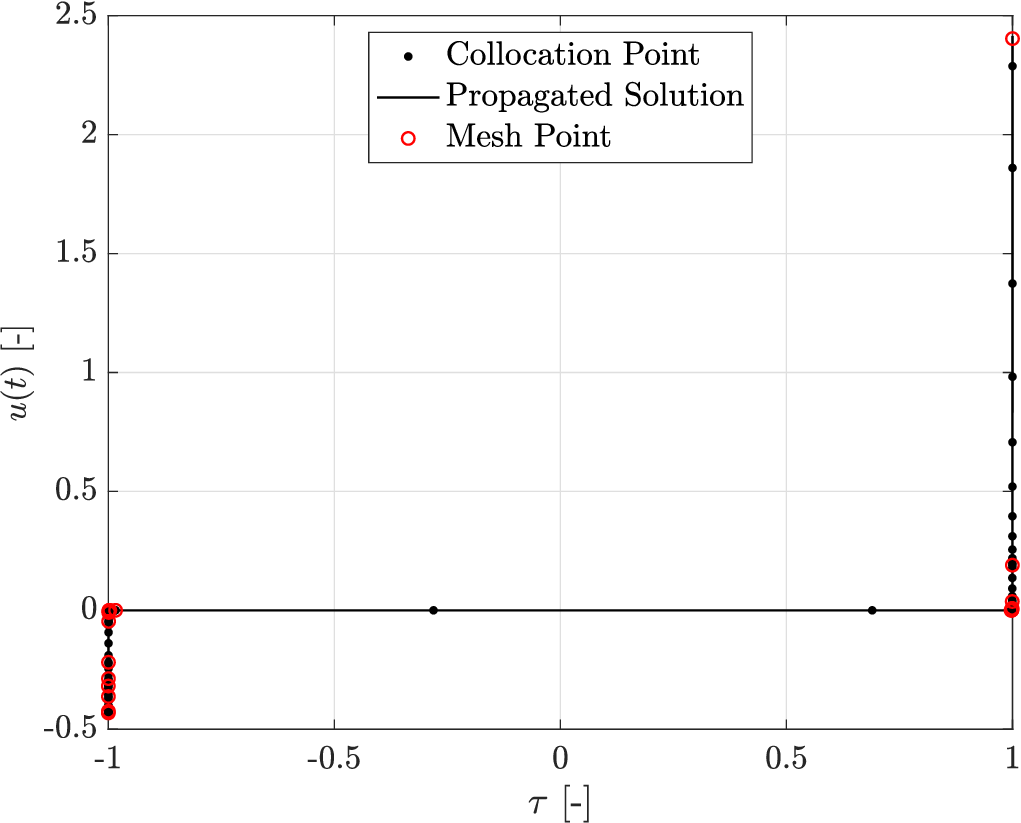}%
        \label{fig:examples-hypersensitive-u}
        }
    \hspace*{\fill}
    \subfloat[Near $\tau=+1$.]{%
        \includegraphics[height=\figureheight]{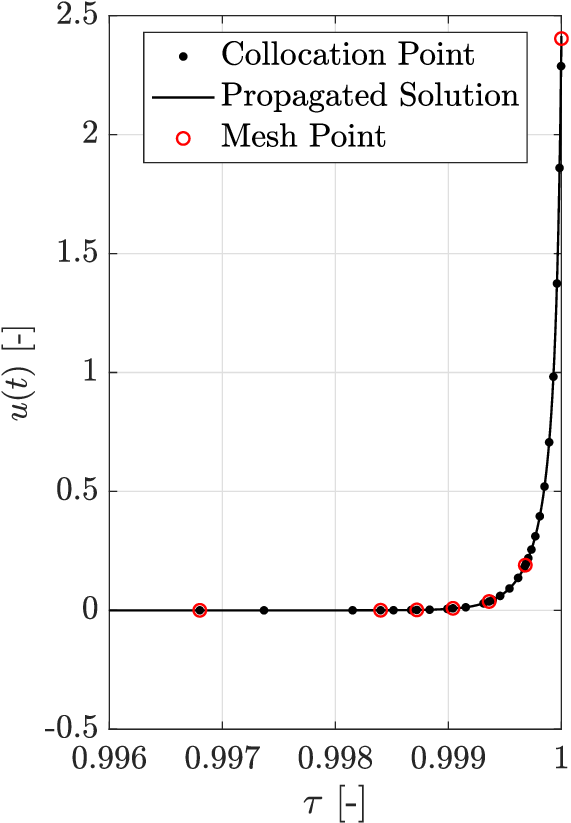}%
        \label{fig:examples-hypersensitive-uend}
        }
    \hspace*{\fill}
    \caption{Optimal control, $u(t)$ vs. $\tau$, for Example 3 using the \HRvar{3}{12}{\ode{89}} method.}
    \label{fig:examples-hypersensitive-control}
\end{figure}

\begin{figure}[!t]
    \centering
    \hspace*{\fill}
    \subfloat[Near $\tau=-1$.]{%
        \includegraphics[height=\figureheight]{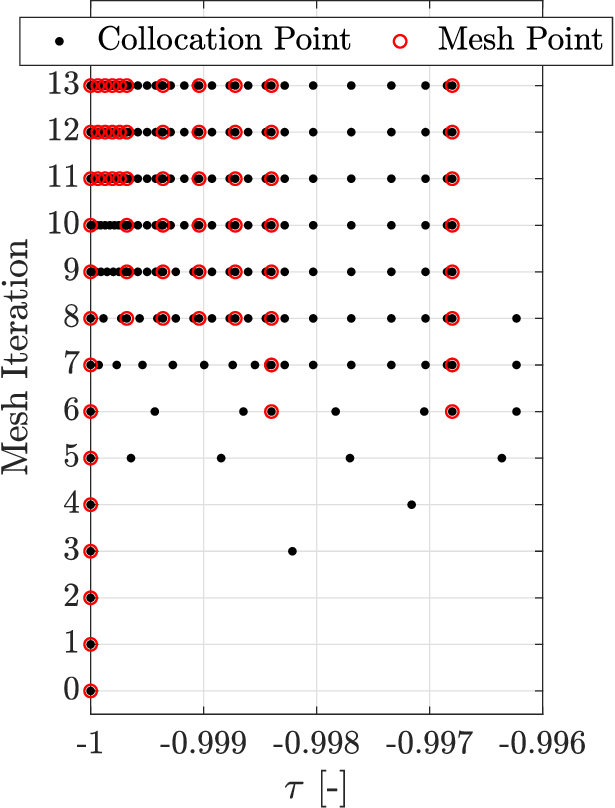}%
        \label{fig:examples-hypersensitive-meshhistorystart}
        }
    \hspace*{\fill}
    \subfloat[Over $\tau\in\ensuremath{[-1,+1]}$.]{%
        \includegraphics[height=\figureheight]{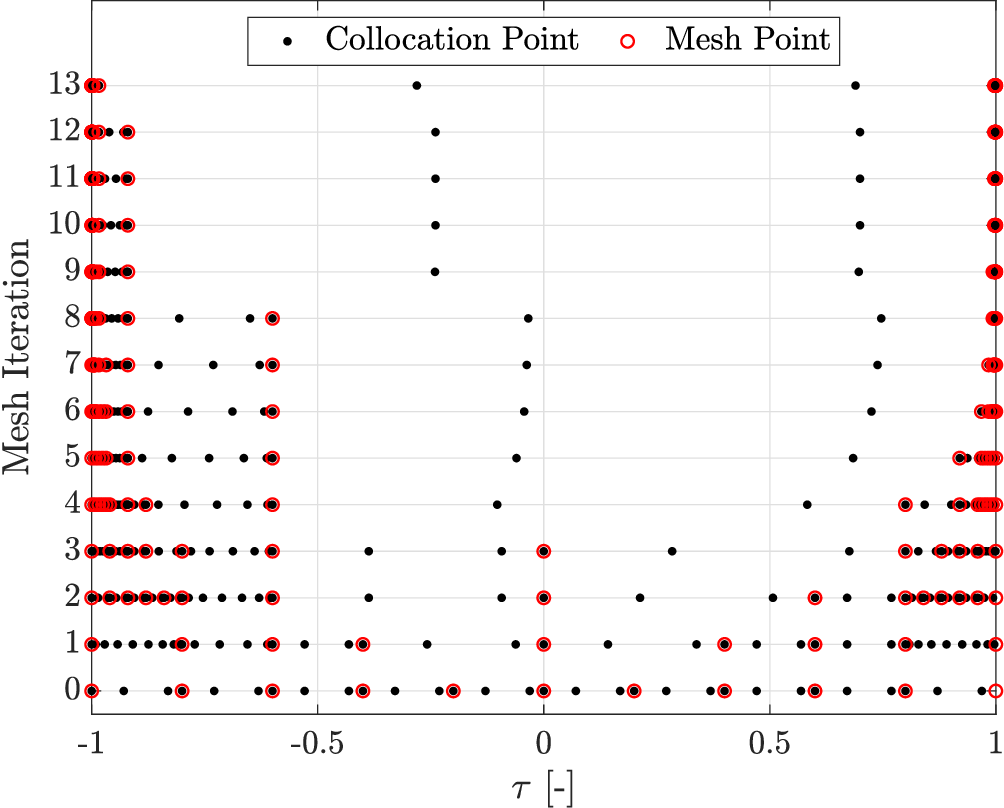}%
        \label{fig:examples-hypersensitive-meshhistory}
        }
    \hspace*{\fill}
    \subfloat[Near $\tau=+1$.]{%
        \includegraphics[height=\figureheight]{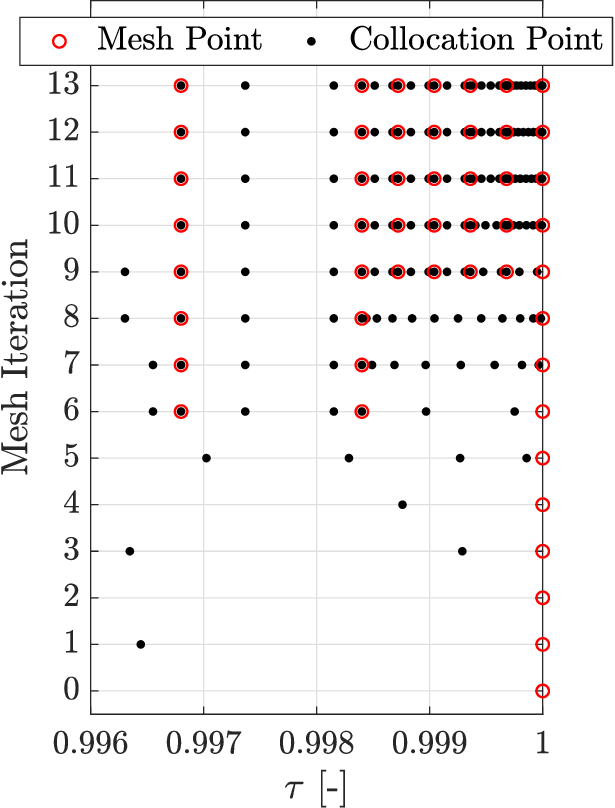}%
        \label{fig:examples-hypersensitive-meshhistoryend}
        }
    \hspace*{\fill}
    \caption{Mesh refinement history for Example 3 using the \HRvar{3}{12}{\ode{89}} method.}
    \label{fig:examples-hypersensitive-mesh_1em6}
\end{figure}

Using the \HRvar{3}{12}{\ode{89}} method, a solution to this optimal control problem on $\tau\in[-1,+1]$  is shown in Figs.~\ref{fig:examples-hypersensitive-state} and~\ref{fig:examples-hypersensitive-control}.
It is noted that Figs.~\ref{fig:examples-hypersensitive-xstart}--\ref{fig:examples-hypersensitive-xend} show the optimal state solution obtained by solving the NLP of Eqs.~\eqref{eq:lgr-obj}--\eqref{eq:lgr-cnt} as well as the solution obtained via explicit simulation using the \MATLAB~ODE solver \ode{89}, whereas Figs.~\ref{fig:examples-hypersensitive-ustart}--\ref{fig:examples-hypersensitive-uend} show the optimal control solution along with the interpolated control used in the explicit simulation.
% ---------------------------------------------------------------------------------------------
% Figure: hyper-sensitive example optimal state components vs tau using phs-(3,12,ode89) method
% ---------------------------------------------------------------------------------------------
% -----------------------------------------------------------------------------------------------
% Figure: hyper-sensitive example optimal control components vs tau using phs-(3,12,ode89) method
% -----------------------------------------------------------------------------------------------
For this solution, an optimal objective value of $\mcl{J}^{*}=1.330806$ is obtained, and the corresponding mesh refinement history is shown in Fig.~\ref{fig:examples-hypersensitive-mesh_1em6}, where the \HRvar{3}{12}{\ode{89}} method takes thirteen (13) iterations to satisfy the desired mesh tolerance.
% -----------------------------------------------------------------------------------------
% Figure: hyper-sensitive example mesh refinement history using the phs-(3,12,ode89) method
% -----------------------------------------------------------------------------------------
It is seen that the optimal state, shown in Fig.~\ref{fig:examples-hypersensitive-state}, attains the aforementioned take-off, cruise, and landing structure, which is reflected in the optimal control solution of Fig.~\ref{fig:examples-hypersensitive-control}.
Because of the rapid initial decay in the take-off region and rapid terminal growth in the landing region, it is expected that the \HR~method increases the density of the mesh in these regions to adequately capture the changes in the solution.
As seen in Fig.~\ref{fig:examples-hypersensitive-meshhistorystart} and~\ref{fig:examples-hypersensitive-meshhistoryend}, the \HR~method is able to identify such regions and concentrate both mesh and collocation points in the regions near $\tau=-1$ and $\tau=+1$.
In the constant portion of the solution in Fig.~\ref{fig:examples-hypersensitive-meshhistory} (that is, not near $\tau=-1$ or $\tau=+1$), the \HR~method consecutively merges mesh intervals in order to reduce the size of the mesh in the interior region of the solution, while retaining an increased mesh density in the initial and terminal regions.

Next, the performance of the various methods and error estimates is compared in a manner similar to that done for Examples 1 and 2 in Sections~\ref{subsec:examples-supersonic} and~\ref{subsec:examples-robotarm}, respectively.
Table~\ref{tab:examples-hypersensitive-results} summarizes the results for Example 3, where only the methods that obtain the smallest final mesh are shown.
% --------------------------------------------------------------------------------
% Table: results for hyper-sensitive example using various mesh refinement methods
\begin{table}[!t]
    \centering
    \caption{Results for Example 3 using various mesh refinement methods.}
    \begin{tabular*}{\textwidth}{@{\extracolsep\fill}clccccccc@{}}
        \toprule\toprule
        \multirow{2}{*}{$N_{\min}$} & \multicolumn{1}{c}{\multirow{2}{*}{Method}} & \multirow{2}{*}{$N$} & \multirow{2}{*}{$K$} & \multirow{2}{*}{$M$} & \multirow{2}{*}{$e_{\max}^{[M]}$} & \multirow{2}{*}{$\tilde{e}_{\max}^{[M]}$} & \multicolumn{2}{c}{Total Time [s]} \\
        &&&&&&& SNOPT & CPU \\
        \midrule\midrule
        \multirow{7}{*}{2} & None\textsuperscript{a}    & 20  & 10 & 0  & 4.787$\times{10}^{-1}$ & 2.796$\times{10}^{2~~}$ & 0.038 & 0.070 \\
        \cmidrule{2-9}     & \LRvar\textsuperscript{a}  & 155 & 32 & 6  & 1.164$\times{10}^{-6}$ & 8.560$\times{10}^{-7}$  & 0.466 & 0.721 \\
                           & \LRHRvar{\ode{45}}         & 172 & 38 & 7  & 2.931$\times{10}^{-7}$ & 2.221$\times{10}^{-7}$  & 0.673 & 1.578 \\
                           & \LRHRvar{\ode{89}}         & 172 & 38 & 7  & 2.931$\times{10}^{-7}$ & 2.221$\times{10}^{-7}$  & 0.669 & 2.640 \\
        \cmidrule{2-9}     & \HRvar{2}{10}{\ode{45}}    & 93  & 24 & 11 & 7.804$\times{10}^{-7}$ & 6.146$\times{10}^{-7}$  & 0.470 & 2.141 \\
                           & \HRvar{2}{10}{\ode{89}}    & 93  & 24 & 11 & 7.798$\times{10}^{-7}$ & 6.146$\times{10}^{-7}$  & 0.480 & 4.442 \\
        \midrule
        \multirow{7}{*}{3} & None\textsuperscript{a}    & 30  & 10 & 0  & 4.196$\times{10}^{-1}$ & 9.570$\times{10}^{1~~}$ & 0.031 & 0.061 \\
        \cmidrule{2-9}     & \LRLvar\textsuperscript{a} & 101 & 25 & 6  & 8.952$\times{10}^{-7}$ & 7.341$\times{10}^{-7}$  & 0.979 & 1.218 \\
                           & \LRLHRvar{\ode{45}}        & 100 & 24 & 6  & 7.444$\times{10}^{-7}$ & 6.149$\times{10}^{-7}$  & 0.327 & 0.925 \\
                           & \LRLHRvar{\ode{89}}        & 100 & 24 & 6  & 7.484$\times{10}^{-7}$ & 6.149$\times{10}^{-7}$  & 0.324 & 1.563 \\
        \cmidrule{2-9}     & \HRvar{3}{8}{\ode{45}}     & 91  & 21 & 12 & 8.320$\times{10}^{-7}$ & 6.825$\times{10}^{-7}$  & 0.606 & 1.765 \\
                           & \HRvar{3}{12}{\ode{89}}    & 87  & 18 & 13 & 7.545$\times{10}^{-7}$ & 6.052$\times{10}^{-7}$  & 0.572 & 5.022 \\
        \midrule
        \multirow{7}{*}{4} & None\textsuperscript{a}    & 40  & 10 & 0  & 3.826$\times{10}^{-1}$ & 4.742$\times{10}^{1~~}$ & 0.058 & 0.088 \\
        \cmidrule{2-9}     & \LRLvar\textsuperscript{a} & 107 & 22 & 6  & 7.925$\times{10}^{-7}$ & 5.506$\times{10}^{-7}$  & 0.400 & 0.637 \\
                           & \LRLHRvar{\ode{45}}        & 109 & 23 & 7  & 8.082$\times{10}^{-7}$ & 6.432$\times{10}^{-7}$  & 0.490 & 1.162 \\
                           & \LRLHRvar{\ode{89}}        & 109 & 23 & 7  & 8.083$\times{10}^{-7}$ & 6.432$\times{10}^{-7}$  & 0.492 & 1.853 \\
        \cmidrule{2-9}     & \HRvar{4}{8}{\ode{45}}     & 94  & 20 & 10 & 8.703$\times{10}^{-7}$ & 7.365$\times{10}^{-7}$  & 0.772 & 1.795 \\
                           & \HRvar{4}{10}{\ode{89}}    & 96  & 20 & 15 & 8.316$\times{10}^{-7}$ & 6.893$\times{10}^{-7}$  & 0.709 & 2.960 \\
        \bottomrule\bottomrule
    \end{tabular*}
    \begin{tablenotes}
        \item[$^{\rm{a}}$] The error estimate $e_{\max}^{[M]}$ of Section~\ref{subsec:mesh-err} is obtained with the \MATLAB~ODE solver \ode{45}.
    \end{tablenotes}
    \label{tab:examples-hypersensitive-results}
\end{table}%
% --------------------------------------------------------------------------------
This example again reiterates trends observed in Sections~\ref{subsec:examples-supersonic} and~\ref{subsec:examples-robotarm}.
First, the \HR~method converges to the smallest mesh for every considered value of $N_{\min}$ when compared against previously developed methods, where a 40\%, 9\%, and 10\% decrease in total number of collocation points and a 25\%, 13\%, and 9\% decrease in total number of mesh intervals are observed for $N_{\min}=[2,3,4]$, respectively.
The \HR~method typically performs more mesh refinement iterations to satisfy the desired mesh tolerance, where the choice of \MATLAB~ODE solver impacts the performance of the \HR~method for $N_{\min}=[3,4]$.
Similar to the previous two examples, the error estimate $e_{\max}^{[M]}$ obtained using the method in Section~\ref{subsec:mesh-err} is larger than the error estimate $\tilde{e}_{\max}^{[M]}$ obtained using the method of Ref.~\cite{PattersonRao2015} on the converged mesh for all results shown in Table~\ref{tab:examples-hypersensitive-results}, where less than an order of magnitude difference between the two error estimates is observed.
On the other hand, the opposite is true on the initial mesh for this example, where a difference of roughly two to three orders of magnitude is observed.
Through this, the error estimate of Ref.~\cite{PattersonRao2015} suggests that there is a much larger discrepancy between solutions obtained via collocation and explicit simulation than actually present, which guides the mesh refinement methods to increase the size of the mesh more than necessary in the appropriate regions.
When using the error estimate of Section~\ref{subsec:mesh-err} in the \HR, \PRHR, \LRHR, and \LRLHR~mesh refinement methods, explicit simulation of the dynamics of Eq.~\eqref{eq:examples-hypersensitive-dyn} in forward time does not fail at any point along the solution on any mesh refinement iteration when using either \ode{45} or \ode{89}.
To summarize for Example 3, the \HR~mesh refinement method again demonstrates an improvement over the \PR, \PRHR, \LR, \LRHR, \LRL, and \LRLHR~mesh refinement methods in terms of obtaining a smaller final mesh that satisfies the desired mesh tolerance.

% ------------------------------------------------------------------------------------------------ %
% --------------------- SUBSECTION: Example 3: Hyper-Sensitive Problem (END) --------------------- %
% ------------------------------------------------------------------------------------------------ %

% ------------------------------------------------------------------------------------------------ %
% ------------------------------------------------------------------------------------------------ %
% ------------------------------------ SECTION: Examples (END) ----------------------------------- %
% ------------------------------------------------------------------------------------------------ %
% ------------------------------------------------------------------------------------------------ %

% ------------------------------------------------------------------------------------------------ %
% ------------------------------------------------------------------------------------------------ %
% ---------------------------------- SECTION: Conclusion (START) --------------------------------- %
% ------------------------------------------------------------------------------------------------ %
% ------------------------------------------------------------------------------------------------ %
\section{Conclusions} \label{sec:conclusion}

An adaptive mesh refinement method for solving optimal control problems using Legendre-Gauss-Radau direct collocation has been described.
In regions of the solution where the desired accuracy tolerance is not satisfied, the method increases the degree of the polynomial approximation within a mesh interval and/or the number of mesh intervals.
The method also allows for mesh size reduction in regions of the solution where the desired accuracy tolerance is satisfied by either merging adjacent mesh intervals or decreasing the degree of the polynomial approximation within a mesh interval.
In addition, the mesh refinement process is driven by a new relative error estimate in the state solution that is based on the differences between the Lagrange polynomial approximation of the state obtained via collocation and forward and backward explicit simulations of the dynamics in each mesh interval.
Using the newly developed error estimate and mesh refinement process, the solutions obtained via collocation and explicit simulation schemes are guaranteed to be in agreement with one other on the final mesh.
The method is demonstrated to three examples from the open literature that highlight various features of the method, where the method is able to handle an active state path constraint and multiple control discontinuities in the solution while the issues associated with failure of the explicit simulation scheme are addressed.
The results of this research show that the size of the final mesh is smaller when compared with previously developed mesh refinement methods.

% ------------------------------------------------------------------------------------------------ %
% ------------------------------------------------------------------------------------------------ %
% ----------------------------------- SECTION: Conclusion (END) ---------------------------------- %
% ------------------------------------------------------------------------------------------------ %
% ------------------------------------------------------------------------------------------------ %

% ------------------------------------------------------------------------------------------------ %
% ------------------------------------------------------------------------------------------------ %
% ------------------------------- SECTION: Acknowledgements (START) ------------------------------ %
% ------------------------------------------------------------------------------------------------ %
% ------------------------------------------------------------------------------------------------ %
\section*{Acknowledgements}

The authors gratefully acknowledge support for this research from the U.S. National Science Foundation under Grant CMMI-2031213 and the NASA Florida Space Grant Consortium under Grant 80NSSC20M0093.

% ------------------------------------------------------------------------------------------------ %
% ------------------------------------------------------------------------------------------------ %
% -------------------------------- SECTION: Acknowledgements (END) ------------------------------- %
% ------------------------------------------------------------------------------------------------ %
% ------------------------------------------------------------------------------------------------ %

% ------------------------------------------------------------------------------------------------ %
% ------------------------------------------------------------------------------------------------ %
% ------------------------------------- Miscellaneous (START) ------------------------------------ %
% ------------------------------------------------------------------------------------------------ %
\renewcommand{\baselinestretch}{1}
\normalsize\normalfont
% ------------------------------------------------------------------------------------------------ %
% -------------------------------------- Miscellaneous (END) ------------------------------------- %
% ------------------------------------------------------------------------------------------------ %
% ------------------------------------------------------------------------------------------------ %

% ------------------------------------------------------------------------------------------------ %
% ------------------------------------------------------------------------------------------------ %
% ---------------------------------- SECTION: References (START) --------------------------------- %
% ------------------------------------------------------------------------------------------------ %
% ------------------------------------------------------------------------------------------------ %

% ------------------------------------------------------------------------------------------------ %
% ------------------------------------------------------------------------------------------------ %
% ----------------------------------- SECTION: References (END) ---------------------------------- %
% ------------------------------------------------------------------------------------------------ %
% ------------------------------------------------------------------------------------------------ %

\end{document}